\documentclass[aps,twocolumn,nofootinbib]{revtex4}%
\usepackage{algorithmic}
\usepackage{graphicx}
\usepackage{amsmath}
\usepackage{amssymb}
\usepackage{subfigure}
\usepackage{epsf}
 \usepackage[dvips]{epsfig}
\newcommand{\real}{\mathbb{R}}
\usepackage{enumerate}
 \usepackage{pst-plot}
 \usepackage{psfrag}
 \usepackage{dcolumn}
\usepackage{bm}

\newcommand{\calS}{{\cal S}}
\newcommand{\calP}{{\cal P}}
\newcommand{\calE}{{\cal E}}

\newcommand{\bY}{{\bf Y}}

\newcommand{\bG}{{\bf G }}
\newcommand{\bA}{{\bf A }}
\newcommand{\bB}{{\bf B }}
\newcommand{\bC}{{\bf C }}
\newcommand{\YY}{{Y}}
\newcommand{\bg}{{\bf g }}
\newcommand{\bc}{{\bf c }}

\newcommand{\bE}{{\bf E}}
\newcommand{\bxi}{{\mbox{\boldmath $\xi$}}}

\newcommand{\bI}{{\bf I}}

\newcommand{\bD}{{\bf D}}
\newcommand{\bS}{{\bf S}}
\newcommand{\calY}{{\cal Y}}
\newcommand{\bt}{{\bf t}}

\newcommand{\eb}{\begin{equation}}
\newcommand{\ee}{\end{equation}}

\newcommand{\bsigma}{\mbox{\boldmath$\sigma$}}

\newcommand{\inte}{\mbox{int}}

\newcommand{\sta}{\mbox{sta}}
\newcommand{\PP}{{P}}

\newtheorem{Lemma}{Lemma}

\newtheorem{Conjecture}{Conjecture}
\newtheorem{Theorem}{Theorem}

\newcommand{\tr}{\mbox{tr}}
\newcounter{algo}
\newenvironment{algo}[1]{\refstepcounter{algo}
\begin{center}
\begin{minipage}{0.9\textwidth}   \hrule\smallskip
\textbf{Algorithm \thealgo: #1}
\par\smallskip\hrule\smallskip\ignorespaces}{\par\smallskip\hrule
\end{minipage}
\end{center}
}{}

\begin{document}

\title{Global Optimal Trajectory in Chaos  and NP-Hardness 
}

\author{Vittorio Latorre}
\altaffiliation{Department of Computer, Control and Management Engineering,
``Sapienza'', University of Rome,  Via Ariosto 25,
00185, Rome Italy}
\email{latorre@dis.uniroma1.it}

\author{David Yang Gao}
\altaffiliation{School of Applied Sciences,
Federation University Australia, Mt Helen, Victoria, 3353, Australia}
              \email{d.gao@federation.edu.au}

\date{\today}

\begin{abstract}
This paper presents a new canonical duality methodology for solving   general nonlinear dynamical systems.
Instead of the conventional iterative methods,
the discretized nonlinear system is first formulated as a global optimization problem via the least squares method.
The canonical duality theory shows that  this nonconvex minimization problem can be solved deterministically in polynomial time
if a global optimality condition is satisfied.  The so-called pseudo-chaos produced by Runge-Kutta type of linear iterations are mainly due to the intrinsic  numerical error accumulations.
Otherwise,  the global optimization problem could be NP-hard and the   nonlinear system can be really chaotic.
     A conjecture is proposed, which reveals the connection between chaos in nonlinear dynamics and NP-hardness in computer science.
     The methodology  and the conjecture are verified by applications to the well-known
      logistic equation, a forced memristive  circuit  and the Lorenz system.
     Computational results show that  the canonical duality theory can be used to identify
     chaotic systems and to obtain  realistic global optimal solutions in
   nonlinear dynamical systems.

\end{abstract}
\keywords{: Chaos \and NP-Hardness \and Nonlinear Dynamics \and Global optimization \and Canonical duality theory}

\maketitle

\section{Problems and Motivation }
We are interested in   a new global optimization  method for solving
the following   nonlinear dynamical system
\begin{equation}\label{eq: diffeq}
(IVP): \;\;
\Bigg\{\begin{array}{rl}
\displaystyle \YY'(t)=&\displaystyle F(t, \YY(t)), \quad t\in [0,T],\vspace{0.5em}\\
\displaystyle \YY(0)=&\displaystyle \YY_0,
\end{array}
\end{equation}
where $\YY:  [0,T]  \rightarrow \real^d$ is a vector-valued  unknown function, $d\ge 1$ is a positive integer,
 $F:  [0,T] \times \real^d\rightarrow \real^d$ is a vector-valued nonlinear operator and $Y_0 \in \real^d$ is a given initial data.
Traditional methods for solving this nonlinear system are based on   Runge-Kutta  iterations.
It is well-known that due to  error accumulation, these traditional iteration methods usually produce  the so-called
``chaotic solutions" for a large class  of nonlinear dynamical systems.
Although  the chaotic behaviors have been   studied extensively during the past 50 years, some fundamental problems are still open, such as an effective description of chaos, rough dependence on initial data,
and the relation between the chaos and computational complexities, etc.

By using  finite difference method and
trapezoidal rule\footnote{Clearly, we can adopt high-order rule for approximation of $F(t,\YY)$ at the $k-1$ step, which will not affect significantly the main results in this paper},   the initial value problem $(IVP)$ in continuous space can be
 discretized in the following nonlinear algebraic system:
\begin{equation}\label{eq: trape}
(\calP_0): \;\;\; \YY_k= \YY_{k-1}+ \frac{\delta}{2}  (F_k + F_{k-1}) ,\quad k=1,\dots, n
\end{equation}
where $\delta > 0$ is a step size, $F_k = F(t_k, \YY_k)$,  and  $n$ is the number of total discretization.
Clearly, this is a nonlinear algebraic system and its unknown $\bY =\{ \YY_k \}  = \{ y_k^i \} \in \real^{d\times n}$ is  a matrix.
Direct methods for solving this nonlinear algebraic system are very difficult.
 If the unknown $\YY_k$ in $F_k$ is replaced by  the iteration $\YY_k = Y_{k-1}+  {\delta} F_{k-1}, $ then $(\calP_0)$ is the well-known modified Euler method.
The popular Runge-Kutta method is also a generalized Euler method.

Rather than  the conventional  iterative approximation from the initial value $Y_0$,  the nonlinear algebraic system  (\ref{eq: trape}) can be precisely
 formulated as a problem of least squares minimization:
\begin{equation}\label{eq:primal}
\displaystyle
(\calP) : \min_{\bY \in \calY_a }  P(\bY)=\frac{1}{2} \sum_{k=1}^{n} \left\| Y_k - Y_{k-1} -\frac{\delta}{2}   (F_k + F_{k-1})   \right\|^2.
\end{equation}
In this paper, $\| * \|$ represents the standard $\ell_2$-norm in $\real^d$, and $\calY_a = \real^{d\times n}$. Clearly, $(\calP)$
 is a global optimization problem.
Due to the nonlinearity of $F(t, Y(t))$, the target function $\PP(\bY)$ is usually nonconvex and the problem could have many local
 and global minimal solutions.
\begin{Lemma}
If $\bY = \{ \YY_k\}$ is a solution of   $(\calP_0)$, then it must be   a  solution of $(\calP)$ and $\PP(\bY)=0$.
\end{Lemma}

\begin{figure*}
  \includegraphics[scale=.55]{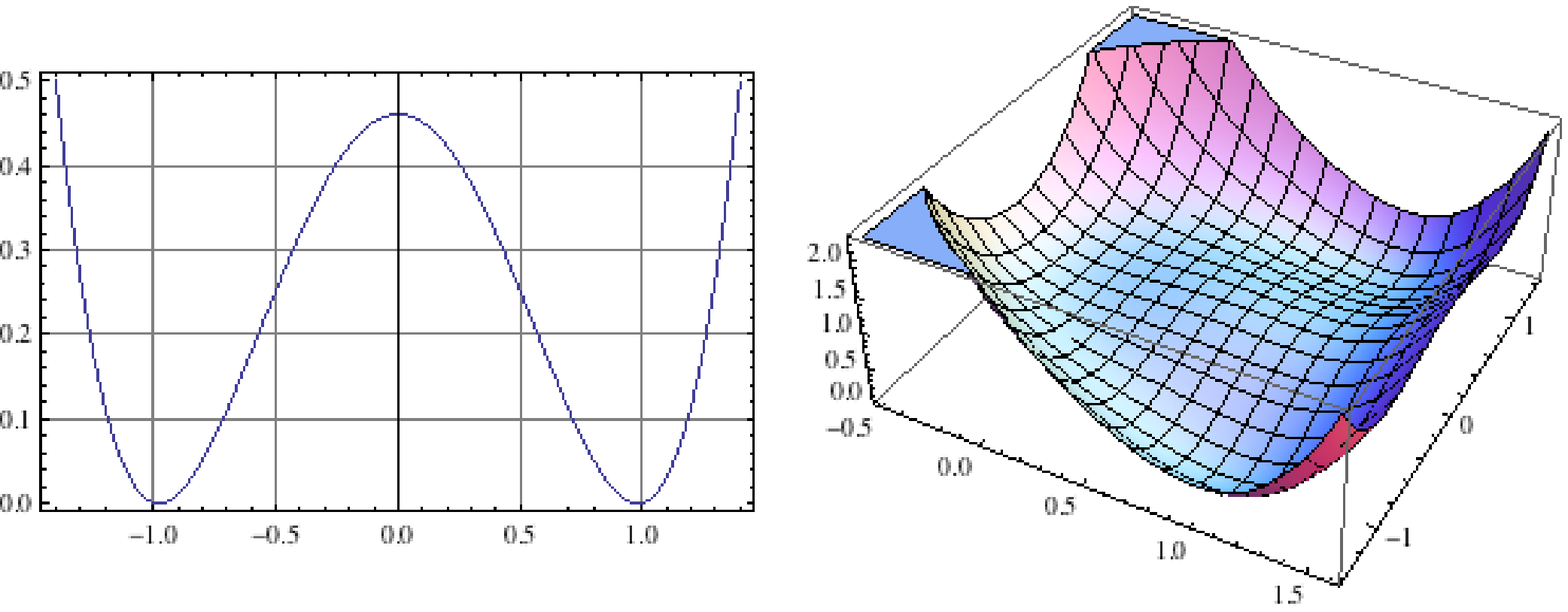}
   $  \mbox{ (a) Graph of $\PP(y_1)$    \hspace{1.5cm} (b) Graph of  $\PP(y_1,y_2)$}$
  \caption{Graphs of $P(Y)$ for logistic equation $y'= r y(1-y)$ with   $r = 4$,  $\delta=0.5 $,  $y_0 = 1.2$.}\label{fig: 1d}
\end{figure*}

This lemma can be proved easily (see \cite{ruan-gao-jiao}), which
 shows that if $\bY = \{ \YY_k\}$ is a solution to  $(\calP_0)$, then it must be a critical point of $\PP(\bY)$.
But not all critical points are solutions to $(\calP_0)$ due to nonconvexity of $\PP(\bY)$.
To see this, let us consider  the  most simple    {\em logistic equation}
\[
y' = r y(1-y), \;\; y(0) = y_0, \;\;  r > 0.
\]
For $n=1$, the target function $P(Y)$ is the so-called {\em double-well function} (see Fig.  \ref{fig: 1d}(a)), which has at most three critical points:
two minimizers  corresponding to the two  algebraic solutions of  the nonlinear  equation:
\eb
y_1 = y_0 + \frac{\delta}{2} r [y_1 (1-y_1) + y_0(1-y_0)],
\ee
and one local maximizer, which is not a solution to the nonlinear equation.
 For $n> 1$, the local minimizers  of $P(Y)$ depend  sensitively on the parameter $r> 0$, the step size $\delta$ and the
 initial data $y_0$ Fig.  \ref{fig: 1d}(b).

The least squares method for solving nonlinear initial value problems was studied by Neuberger and Renka \cite{neu-ren} in 2004.
By using a method of steepest descent with a Sobolev gradient, their  numerical method  for
 Lorenz equations  produces a simple smooth curve that terminates at a stationary point rather than producing the 'butterfly' attractor and they claimed: {\em
As  far  as  we  know this  is  the  first  method  that  produces  a  non-chaotic orbit.}

Generally speaking, if the nonlinear algebraic equation (\ref{eq: trape}) has $K$ solutions for each $Y_k$,   the nonconvex target function
$P(Y)$ could  have about  $K^n$   minimizers.
Conventional optimization methods   such as   Newton or quasi-Newton type of techniques  are   local approaches,
which can only assure the convergence to one of local  solutions.
Similar to the Runge-Kutta iterative method, these  local solutions depend  sensitively on initial data and numerical errors.
It was discovered by Gao and Ogden in \cite{gao-qjmam08} that if a nonlinear ODE system is subjected to
an alternative external force field, both the local and global minimum  solutions could be nonsmooth and can't be captured by  Newton-Type direct methods.
Due to the lack of global optimality criteria, most nonconvex minimization problems are considered to be NP-hard (Nondeterministic method Polynomial-time Hard).
Therefore, from  the  point view  of global optimization and computational complexity theory,
 we now  easily understand that  Neuberger and Renka's smooth curve is only a local stationary solution to the nonconvex minimization problem (\ref{eq:primal}).
How to  solve a general nonconvex minimization problem has been a fundamentally challenging task not only in global optimization and computer science,
but also in multidisciplinary fields of complex systems.


Canonical duality theory is a potentially powerful methodological theory developed from  nonconvex analysis/mechanics  and global optimization \cite{gbook,gao09}.
This theory can be used not only to model complex systems within a unified framework, but also for solving a large class of nonconvex/nonsmooth/discrete
optimization problems in multidisciplinary fields of computational biology, ecology \cite{ruan-gao-ima2014},
engineering sciences \cite{gao09,wang12}, and recently in network communications \cite{latorreSNL,ruan-gao-pe2014},  nonconvex constrained optimization \cite{lg15,ls14}, and radial basis neural networks \cite{lg14}, etc. Comprehensive reviews on the canonical duality theory are given in \cite{gao-cace09,gao-ruan-la14}.
The main goal of this paper is to apply the canonical duality theory for solving the nonconvex minimization problem (\ref{eq:primal}).
The next section provides detailed information on   how to reformulate $(\calP)$ as a  concave maximization  dual problem and
what is the global optimality condition.
Numerical method and perturbation technique are discussed in Section 3.
Based on the canonical duality theory, a conjecture is proposed which reveals the connection between chaos in nonlinear dynamics and NP-hardness in computer science.
Applications are illustrated in Section 4.
The methodology and conjecture are tested by  three well-known systems. The first one is one-D logistic equation, it has a unique global optimal trajectory.
The second one  is a two-D forced memristive circuit.  This system possesses chaotic behavior by conventional iterative method   \cite{itoh11}. But by the canonical duality theory, this system has
global optimal solutions  for different initial data and parameters.
The final example is the classical 3-D Lorenz system with the butterfly chaotic attractor. The canonical duality theory shows reason
why  this well-known chaotic system is NP-hard.
Conclusions are presented in the last section.


\section{Canonical Duality Theory and Global Optimality Condition}
The canonical duality theory is composed mainly of (1) a canonical transformation; (2) a complementary-dual principle;
(3) a triality theory.
To simplify the exposition,  we assume that the nonlinear operator  $ F (t , Y )$  is a vector-valued  quadratic function of $Y$:
\[
 F(t ,Y ) =\YY^T  \bA  \YY + \bD\YY +  \bg(t ) ,
\]
where $\bg(t) = \{ g_i(t) \} \in \real^d$ is a given vector-valued function,
$\bA = \{ \bA^i_{jl}\}  \in\real^{d\times d\times d}  $ and $\bD= \{D_{ij} \}  \in\real^{d\times d} $ are given generalized  matrices, respectively.
 Thus, the target function $\PP(\bY)$ is a fourth-order polynomial of $\bY \in \real^{d\times n}$:
  \begin{equation}\label{eq:primal1}
 \PP(\bY) =
 \frac{1}{2}\sum_{k=1}^{n} \left\| \Lambda(\YY_k) + \bB (\YY_k)  + \bc_k    \right\|^2 ,
\end{equation}
where $\bB(\YY_k) =\frac{\delta}{2}   \bD (\YY_{k}+ \YY_{k-1})  -  \YY_k + \YY_{k-1} $, $\bc_k = \frac{\delta}{2} (\bg_k + \bg_{k-1})$, and
\eb
\Lambda(\YY_k) =   \frac{\delta}{2}  ( \YY_{k}^T \bA  \YY_{k}+ \YY_{k-1}^T \bA  \YY_{k-1}) .
\ee

The key idea of the canonical dual transformation is to choose a certain geometrically admissible  non-linear measure such that the high-order
nonconvex target function can be written in a  canonical function in dual space.
By the definition used in the canonical duality theory, a real-valued function $V:\calE_a \rightarrow \real$ is said to be a canonical function on its domain $\calE_a$
if the canonical dual mapping $\bsigma = \nabla  V(\bxi) : \calE_a \rightarrow \calE^*_a$ is one-to-one and onto.
For the fourth-order polynomial $\PP(\bY)$, such a geometrical operator can be chosen as
\eb
\bE = \{ \bxi_k\} = \{\xi^i_{k}\} =\Lambda(\bY) = \{ \Lambda(\YY_k)\}:  \calY_a  \rightarrow{\cal E}_a\subseteq\real^{d\times n}.
\ee
Thus, a bi-quadratic (canonical) function   $ \real^{d\times n} \times \calE_a \rightarrow \real$ can be defined by
\eb
\Phi(\bY, \bE) = \frac{1}{2}\sum_{k=1}^{n} \left\| \bxi_k  + \bB (\YY_k)  + \bc_k    \right\|^2 .
\ee
In terms of the geometrical operator $\bE = \Lambda (\bY)$, the primal problem $(\calP)$ can be equivalently written in the canonical form
\begin{equation}\label{eq:primtra}
(\calP): \;\;\; \min \{ \Pi(\bY)= \Phi(\bY, \Lambda(\bY)) \; | \;\; \forall \bY \in \calY_a \} .
\end{equation}

 By the fact that for any given $\bY \in \calY_a$, the canonical function $\Phi(\bY, \bE)$ is convex in $\bE$, the canonical dual variable
 $\bS \in \calE^*_a \subset \real^{d\times n}$ can be uniquely defined by
 \eb
 \begin{array}{rl}
 \bS& \displaystyle=\{ \bsigma_k\} = \{ \sigma^i_k\}  =
  \nabla_{\bE}  \Phi(\bY, \bE) \vspace{0.5em}\\
  &= 
  \bE + \bB(\bY) + \bC = \{ \bxi_k + \bB(\YY_k) + \bc_k\},
  \end{array}
 \ee
 where $\calE^*_a$ is the range of $\nabla_{\bE} \Phi(\bY, \bE)$.
For a given $\bY \in \calY_a$, the conjugate function $\Phi^*(\bY, \bS)$ can be defined by the (partial)  Legendre transformation
 \begin{eqnarray}
 \Phi^*(\bY,\bS) &=&  \{ \tr(\bE^T \bS) -  \Phi(\bY, \bE) \; | \;\; \bS = \nabla_{\bE} \Phi(\bY, \bE)  \} \nonumber \\
& =&\ \sum_{k=1}^{n}  \left[ \frac{1}{2} \| \bsigma_{k}\|^2  - \bsigma^T_k \bB(\YY_k)  - \bc_k^T \bsigma_k \right] ,
\end{eqnarray}
 where $\tr = $trace. Clearly, for any given $\bY \in \calY_a$, the conjugate $\Phi^*(\bY, \bS)$ is a canonical function of $\bS$ and we have the following
 canonical duality relations:
 \begin{eqnarray}
& &  \bS = \nabla_{\bE} \Phi(\bY, \bE) \;\;\Leftrightarrow \;\;  \bE = \nabla_{\bS} \Phi^*(\bY, \bS) \;\; \Leftrightarrow  \nonumber \\
&  &  \Phi(\bY, \bE) + \Phi^*(\bY, \bS) = \tr (\bE^T \bS) \;\; \forall \bY \in \calY_a. \label{eq-fye}
 \end{eqnarray}
 The equation (\ref{eq-fye}) is the so-called Fenchel-Young equality, by which, the canonical function $\Phi(\bY, \Lambda(\bY))$ can be written in
 Gao-Strang's total complementary function:
\begin{eqnarray}\label{eq: xi}
    \Xi(\bY,\bS) &=& \tr(\Lambda(\bY)^T \bS) -  \Phi^*(\bY, \bS)   \} \nonumber \\
& =&  \sum_{k=1}^n  \Bigg[ \frac{1}{2} \YY^T_k \bG(\bsigma_k) \YY_k - \YY^T_k \bt(\bsigma_k)\nonumber\\
   &&- \frac{1}{2}\|\bsigma_k \|^2  + \bc^T_k \bsigma_k  \Bigg] + \psi(\bsigma_1), \;\;\; \;\;
\end{eqnarray}
where $\bG  \in \real^{d\times d}$, $\bt\in \real^d$, and $\psi(\bsigma_1) \in \real$ are  defined by
\begin{eqnarray}
  \bG(\bsigma_k)&=&  \left\{  \sum_{i=1}^d \bA^i (\sigma^i_{k}+\sigma^i_{k+1})  \right\} \\
  && \;\; \forall  \;\;  k=1,\dots, n-1 \label{eq:gsig}\\
  \bG (\bsigma_n) &=& \left\{ \sum_{i=1}^d \bA^i \sigma^{i}_{n} \right\}\\
  \bt (\bsigma_k) &=& \left\{ \sigma^{i}_{ k}-\sigma^i_{ k+1}   -  \frac{\delta}{2}\sum_{j=1}^{d} D_{ji}(\sigma^{j}_{ k+1}+\sigma^{j}_{k}) \right\}\nonumber\\
  && \;\;  \forall   i=1,\dots, d, \;\;  \forall   k=1,\dots, n-1 \\
   \bt(\bsigma_n) &=& \left\{  \sigma^{i}_{n}   -   \frac{\delta}{2}\sum_{j=1}^d  D_{ji}\sigma^{j}_{n}\right\}  \;\;  \forall   i=1,\dots, d\\
   \psi(\bsigma_1) &= &   \YY_{0}^T\bsigma_{1} +  \frac{\delta}{2}(\YY_0^T \bA \YY_0 + \bD  \YY_0)\bsigma_{1}.
\end{eqnarray}

 It is easy to see that $\Xi:\calY_a \times \calE^*_a \rightarrow \real$ is a bi-quadratic function, by which,
  the canonical dual function can be obtained by
\begin{eqnarray}
\Pi^d(\bS) &=& \left\{ \Xi(\bY, \bS) | \;\; \nabla_{\bY} \Xi(\bY, \bS) = 0 \;\;\forall \bS \in \calE^*_a \right\} \nonumber \\
&=&  \sum_{k=1}^n  \Bigg[  - \frac{1}{2} \bt(\bsigma_k)^T \bG(\bsigma_k)^{-1} \bt(\bsigma_k )\nonumber\\
&&- \frac{1}{2}\|\bsigma_k \|^2+  \bc^T_k \bsigma_k   \Bigg] +\psi(\bsigma_1). \;\;\;\;\; \label{eq: dual}
\end{eqnarray}
 Let $\calS_a = \{ \bS \in \calE^*_a| \;\; \det \bG(\bsigma_k) \neq 0 \;\;\forall k = 1, \dots, n\}$ be an admissible  dual space, on which,  the canonical dual function  $\Pi^d(\bS)$ is well-defined.
\begin{Theorem}
\label{th: tcp}
$(${\bf Complementarity-Dual Principle \cite{gbook}}$)$
The function $\Pi^d(\bS)$ is canonically dual to $\Pi(\bY)$ in the sense that if $\bar\bS$ is a critical point of $\Pi^d(\bS)$ then
\begin{equation}\label{eq: xinsig}
\bar{\bY} ({\bar \bS}) = \{ \bG^{-1}(\bar\bsigma_k) \bt(\bar\bsigma_k) \}
\end{equation}
is a critical solution to $(\calP)$ and
$$
\Pi(\bar{\bY})=\Xi(\bar{\bY},\bar\bS)= \Pi^d(\bar\bS).
$$
Conversely, if $\bar{\bY}$ is a solution to  the nonlinear algebraic system $(\calP_0)$,
it must be in the form (\ref{eq: xinsig}) for certain critical point   $\bar\bS$ of $\Pi^d(\bS)$.
\end{Theorem}

This theorem shows that there is no duality gap between the nonconvex target function $\Pi(\bY)$ and its canonical dual $\Pi^d(\bS)$, and the local solutions to the
nonlinear algebraic system $(\calP_0)$ can be analytically represented in the canonical dual solution.
Note that the dual feasible set  $\calS_a$ is  not convex, then
 in  order to find both local and global minimum solutions,  we need to introduce the following subsets of $\calE^*_a$:
\begin{eqnarray}
\calS^+_a &=& \{\bS = \{ \bsigma_k \} \in \calE^*_a|\; \bG(\bsigma_k) \succeq  0 \;\; \forall k = 1, \dots, n\}, \nonumber\\
 \calS^-_a & =&  \{\bS = \{ \bsigma_k \} \in \calE^*_a|\; \bG(\bsigma_k)  \prec   0 \;\; \forall k = 1, \dots, n \}.\nonumber
\end{eqnarray}
By the fact that the block matrix $\bG(\bsigma_k)$ in $\calS^+_a$ is allowed to be singular ($ \det \bG(\bsigma_k) = 0 $),
the term $\bG(\bsigma_k)^{-1}$ in $\Pi^d(\bS)$ should be understood as the Moore-Penrose pseudoinverse \cite{gao-cace09}.
\begin{Theorem}\label{th: triality} $(${\bf Triality Theory \cite{gbook}}$)$
Let  $(\bar{\bY},\bar\bS)$ be an isolated  critical point of $\Xi(\bY,\bS)$.
The following  three extremality conditions hold:
\begin{enumerate}
\item {\bf Global Optimum}: The critical point  $\bar{\bY}$ is a
 global minimizer of $\Pi(\bY)$ if and only if $\bar\bS\in\calS_a^+$. In this case, we have
\begin{equation}\label{eq-tris}
\,\,\,\,\,\,\,\,\,\,\,\,\,\,\,\,\min_{\bY\in \calY_a}\Pi(\bY)=\Pi(\bar{\bY})= \Pi^d(\bar\bS)=\max_{\bS\in\calS_a^+}\Pi^d(\bS).
\end{equation}
\item {\bf Local Maximum}: If $\bar\bS \in\calS_a^-$,  then $\bar\bS $ is a local maximizer of $\Pi^d(\bS )$ on its neighborhood $\calS_o\subset \calS_a^-$ if and only if $\bar{\bY}$ is a local maximizer of $\Pi(\bY)$ on its neighborhood ${\calY}_o\in\calY_a$. In this case, we have
\eb \label{eq-trima}
\,\,\,\,\,\,\,\,\,\,\,\,\,\,\,\,\,\max_{\bY\in {\calY}_o}\Pi(\bY)=\Pi(\bar{\bY})= \Pi^d(\bar\bS)=\max_{\bS\in\calS_o}\Pi^d(\bS).
\ee
\item {\bf Local Minimum}: If $\bar\bS \in\calS_a^-$,  then $\bar\bS $ is a local minimizer of $\Pi^d(\bS )$ on its neighborhood $\calS_o\subset \calS_a^-$ if and only if $\bar{\bY}$ is a local minimizer of $\Pi(\bY)$ on its neighborhood ${\calY}_o\in\calY_a$. In this case, we have
\eb \label{eq-trimi}
\,\,\,\,\,\,\,\,\,\,\,\,\,\,\,\,\,\min_{\bY\in {\calY}_o}\Pi(\bY)=\Pi(\bar{\bY})= \Pi^d(\bar\bS)=\min_{\bS\in\calS_o}\Pi^d(\bS).
\ee
\end{enumerate}
\end{Theorem}

The first statement (\ref{eq-tris}) is  called {\em canonical min-max duality}. Its
  was developed from  Gao and Strang's work  in 1989 \cite{gao-strang89}. This duality
 can be used to identify global minimizer of the nonconvex problem $(\calP)$.
According this statement, the nonconvex problem $(\calP)$  is equivalent to
the following canonical dual problem, denoted by $(\calP^d)$:
\eb
(\calP^d): \;\;\; \max \sta \{ \Pi^d(\bS) | \; \bS \in \calS_a^+  \} , \label{eq-cdmax}
\ee
i.e. to  find the global maximizer among all stationary points of $\Pi^d(\bS)$ on $\calS^+_a$.
By the fact that the canonical dual   $\Pi^d(\bS)$  is a strictly concave function  over a close convex domain $\calS^+_a$,
this canonical dual problem  can be solved easily by well-developed convex analysis and optimization techniques
if $\calS^+_a \neq \emptyset$.
The solution obtained  by this statement provides  a global optimal trajectory for the nonlinear system $(IVP)$.

The second statement (\ref{eq-trima}) is the {\em canonical double-max duality} and the third one
  (\ref{eq-trimi}) is the  {\em canonical double-min duality}.
  These two statements can be used to identify the two special local extremum solutions $\bY^\sharp$ and $\bY^\flat$ such that
 $\Pi(\bY^\sharp) \ge \Pi(\bar \bY)$ for all stationary solutions ${\bar \bY}$ and $\Pi(\bY^\flat) \ge \Pi({\bar \bY})$ for all
  local minimizers ${\bar \bY}$.
  By the fact that $\bY^\sharp$ is not a solution to $(\calP_0)$, the local minimum solution $\bY^\flat$
   determined by  (\ref{eq-trimi}) could be useful for
  understanding the most unstable trajectory in chaos.

\section{Numerical Methodology and Criteria for Chaos and NP-hardness}

This section has two goals: 1) to propose a methodology for solving the canonical dual problem $(\calP^d)$; 2) to find criteria for
identify chaos and NP-hard problems.

 Theorem \ref{th: triality} shows that there are different approaches for solving  the nonconvex primal problem $(\calP)$.
One  is to directly solve the dual problem  $(\calP^d)$, but this approach  has two main  disadvantages:
\begin{itemize}
\item It is necessary to calculate the inverses of matrices $\bG(\bsigma_k)$ for $k=1,\dots,n$ every time  the target function is evaluated,
 and such operation could be necessary several times per iteration;
\item To compute the inverse matrix  can  be time expensive or generate errors in the case even one of the  $\bG(\bsigma_k)$ for $k=1,\dots,n$ is  either  ill-conditioned or   not full rank.
\end{itemize}
For these reasons we propose a canonical primal-dual method for solving the challenging nonlinear system $(IVP)$.

From Theorem \ref{th: tcp} it is clear that there is a one to one correspondence among  the stationary points of  $\Pi(\bY)$, $\Pi^d(\bS)$ and
$\Xi(\bY,\bS)$.
Therefore it is possible to find a global optimum for Problem (\ref{eq:primal}) by solving the following saddle-point problem ($(SP)$ in short):
\begin{equation}\label{eq: xigen}
(SP): \;\;\; \min_{\bY\in\calY_a }\max_{\bS\in\calS_a^+}\Xi( \bY,\bS)= \max_{\bS\in\calS_a^+} \min_{\bY\in\calY_a }\Xi( \bY,\bS).
\end{equation}

Methods for solving  saddle point problems have been extensively  studied in literature \cite{saddle}.
One of the most used approaches is to reformulate the  (\ref{eq: xigen}) as a variational inequality problem on a
closed convex set (cone): 
\begin{equation}\label{eq: VIgen}
\Gamma (\bY,\bS)=\left(
\begin{array}{c}
\nabla_{\bY} \Xi(\bY,\bS) \\
-\nabla_{\bS}  \Xi(\bY, \bS) \\
\end{array}\right) = 0  \;\; \forall\;  \bS \in \calS^+_a .
 \ee
By the fact that the operator $\Gamma$ is monotone  and the feasible set  $\calS^+_a$ is convex,
this problem is equivalent to a convex minimization  and the solution of (\ref{eq: VIgen}) can be obtained easily.

In the case that the canonical dual solutions are located on the boundary of   $\calS_a^+$, perturbation methods can be suggested \cite{gao-ruan-jogo10,ruan-gao-pe}.
   The simplest perturbation problem   is the following:
\begin{equation}\label{eq: perprim}
(\calP_\rho) : \;\; \Pi_\rho(\bY)= \Pi(\bY) + \frac{1}{2} \rho \sum_{k=1}^{n}  \|\YY_k \|^2,
\end{equation}
where $\rho > 0$ is a perturbation parameter.
In this problem, the  canonical dual feasible set $\calS^+_a$ should be replaced by
\eb
\begin{array}{rl}
\calS^+_\rho =& \{ \bS \in\calE^*_a | \;\; \bG_\rho(\bsigma_k )=  \vspace{0.5em}\\
&\bG(\bsigma_k)+  \rho \bI_d \succeq 0 \;\; \forall k = 1, \dots, n\},
\end{array}
\ee
where $\bI_d $ is an identity matrix in $\real^{d\times d}$.
It is easy to see  that as the perturbation
parameter $\rho$ is getting larger and larger, the matrices  $\{ \bG_\rho(\bsigma_k)\} $ become  more and more diagonally dominated and positive definite
and the primal problem $\Pi_\rho(\bY)$ becomes convex.
 On the other hand as  $\rho$ is approaching to zero, the perturbed problem and its dual  are approaching  to their  original formulations.
 Such behavior can be exploited in an optimization method to find a good approximation of the original problem $\Pi(\bY)$. As a matter of facts it could be possible to start the optimization with a large enough value of $\rho$ and then lower this value during the iterations in order to reach a good enough approximation   solution of $\Pi(\bY)$.

 The triality theory 
and the associated  canonical primal-dual method   have been applied successfully for solving many challenging problems
 in computational sciences (see \cite{ruan-gao-ima2014,ruan-gao-pe,zhang-gao-jtb}).
 For the nonlinear dynamical system $(IVP)$, if its  discretized nonlinear algebraic problem $(\calP_0)$ has
 multiple  solutions for certain  $\YY_k$, then the primal problem $(\calP)$ could have  multiple global  minimizers.
Conventional iteration methods for solving $(\calP_0)$
 could produce chaos since these methods are sensitive to the initial data and numerical error accumulation.
However, by the triality  theory  we know that 
 the nonconvex minimization problem $(\calP)$ is not NP-hard if
its canonical dual problem $(\calP^d)$ is solvable, i.e. it has solutions in $\calS^+_a$. These solutions are usually located on the
boundary of $\calS^+_a$, which can be obtained deterministically by perturbation method \cite{ruan-gao-pe} if   $\inte \calS^+_a \neq \emptyset$,   i.e.   $\calS^+_a$ has  interior.
Dually, if $\inte \calS^+_a= \emptyset$, the canonical dual problem $(\calP^d)$  is not solvable. In this case,
   the complementary-dual principle (Theorem 1)
shows that the primal problem $(\calP)$ is equivalent to the following minimum stationary problem:
\eb  \label{cd-s}
(\calP^d_s): \;\; \min \sta \{ \Pi^d(\bS ) | \;\; \bS \in \calS_a \},
\ee
i.e. to find the global minimizer among all stationary points of $\Pi^d(\bS)$ on $\calS_a$.
This is a nonconvex minimization problem over a nonconvex feasible space, which could be really NP-hard.
Therefore, a conjecture was proposed in \cite{gao07}.
\begin{Conjecture}[NP-hard Criterion]\label{np-hard} 
The problem $(\calP)$ is NP-Hard if its canonical dual $(\calP^d)$ is not solvable.
\end{Conjecture}


This  conjecture 
is important for understanding chaos in nonlinear dynamical systems.
By looking at the definitions of $\calS_a$ and $\calS_a^+$,  we know  that if even one of the matrices $\bG(\bsigma_k)$ for $k=1,\dots, n$  is  indefinite for any given  $\bS = \{\bsigma_k\} \in \calE^*_a$, then $\inte \calS^+_a = \emptyset$.
Also, the existence of the canonical dual solution to $(\calP^d)$ depends  on the system  inputs   $\bg(t)$ (see \cite{ruan-gao-pe}).
 Such empirical behaviors for Problem $(\calP)$  and its dual $(\calP^d)$  give rise to a  new conjecture:
\begin{Conjecture}[Chaos and NP-Hardness]\label{con: c=np}
The nonlinear system  $(\calP_0)$  has   chaotic solutions
 if and only if the optimization problem $(\calP)$  is NP-Hard, i.e.  its
canonical dual $(\calP^d)$ is not solvable. 
\end{Conjecture}

This conjecture reveals  the connection between chaos in nonlinear dynamics and NP-hardness in computer science.
Based on  this conjecture,  the triality  theory
 can be used to identify chaotic behaviors of a nonlinear dynamical system  in the following ways.

 1) If  $(\calP^d)$ has a unique  solution $\bar\bS \in \calS^+_a$, 
  the nonlinear system $(IVP)$  is stable and the conventional iteration methods for solving $(IVP)$ shouldn't  produce chaos.

  2) If   $(\calP^d)$ has multiple solutions on the boundary of $\calS^+_a$ and $\inte \calS^+_a \neq \emptyset$, 
   the nonlinear system $(IVP)$ can be considered as {\em deterministically stable} since its global optimal trajectories can be
  obtained deterministically by the canonical duality theory.
  But,  the conventional iteration methods may produce chaotic ``solutions". This  type of chaotic behavior  can be called {\em pseudo-chaos} which is intrinsic to the numerical  error-accumulation of these iteration methods.

  3)  If $(\calP^d)$ is not solvable, the problem $(\calP)$ is NP-hard and the  nonlinear system $(IVP)$ is chaotic.

The methodology, conjecture and results presented in this section can  be  numerically verified by  examples   in  the next section.

\section{Numerical Experience}
This section presents applications of  the proposed theory and methodology via  three examples:
\begin{enumerate}
\item One-D Logistic equation with unique solution in  $\calS^+_a$;
\item Two-D forced memristive circuit  with   $\inte \calS^+_a \neq \emptyset$;
\item Three-D  Lorenz system  with  $\inte \calS^+_a = \emptyset $.
\end{enumerate}
As a benchmark we use the standard ode45 routine from matlab, which is based on the Runge-Kutta (R-K) method. For every single instance, we first run the test with ode45 and then use the result as starting  point of iteration for the Canonical Duality (CD) methodology.
The quality of the results is ascertained by plugging back the numerical results obtained by the two methodologies into the primal function (\ref{eq:primal}). The one  with the lower value of the target  function is considered to be the trajectory closer to the true solution.

\subsection{Logistic Equation}
\begin{figure}
\subfigure[Trajectory by CD methodology.]{\includegraphics[scale=.18]{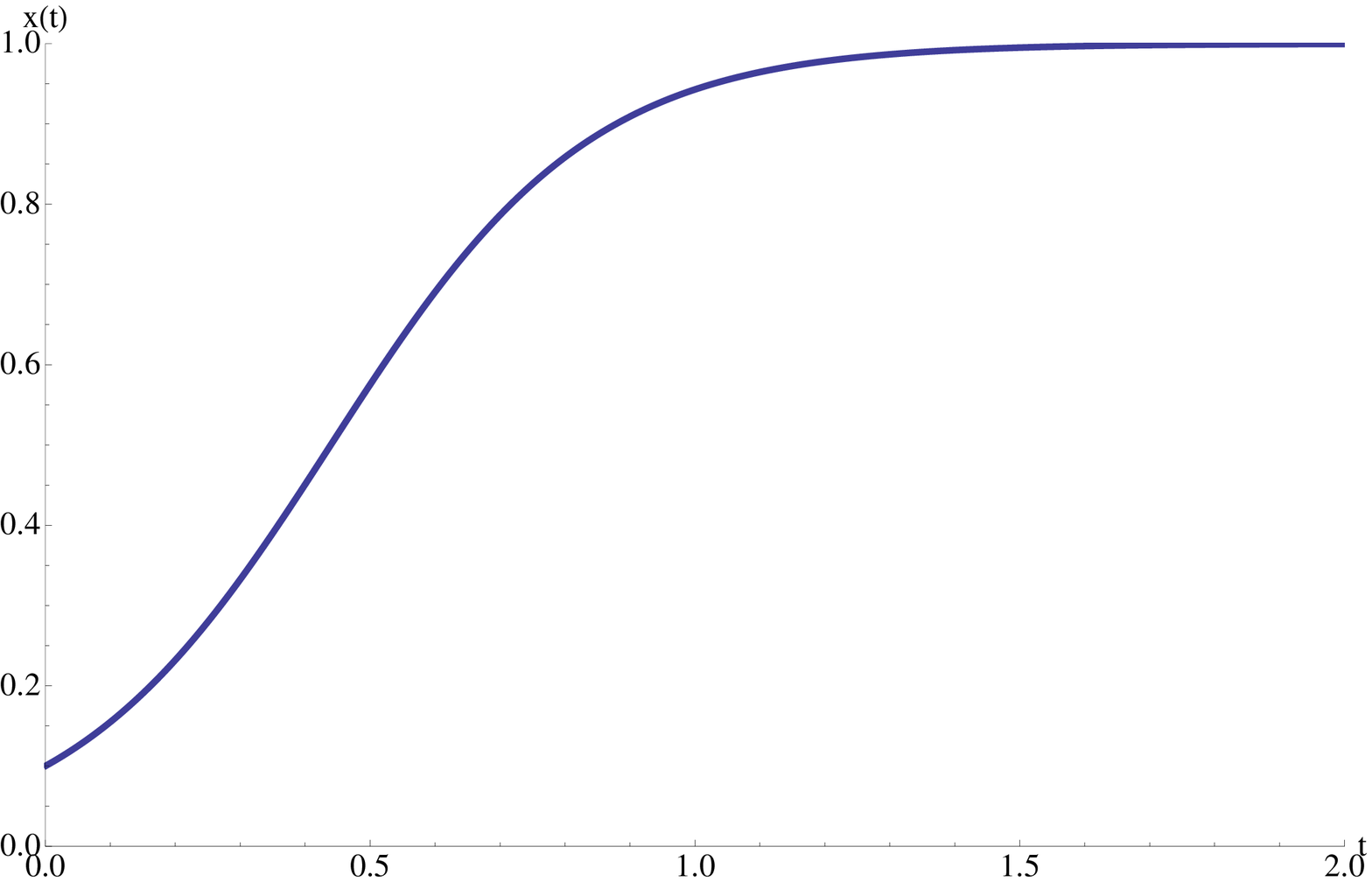}\label{fig: 1dcan}}
\hspace{.5cm} \subfigure[Trajectory  by R-K (ode45).]
{\includegraphics[scale=.18]{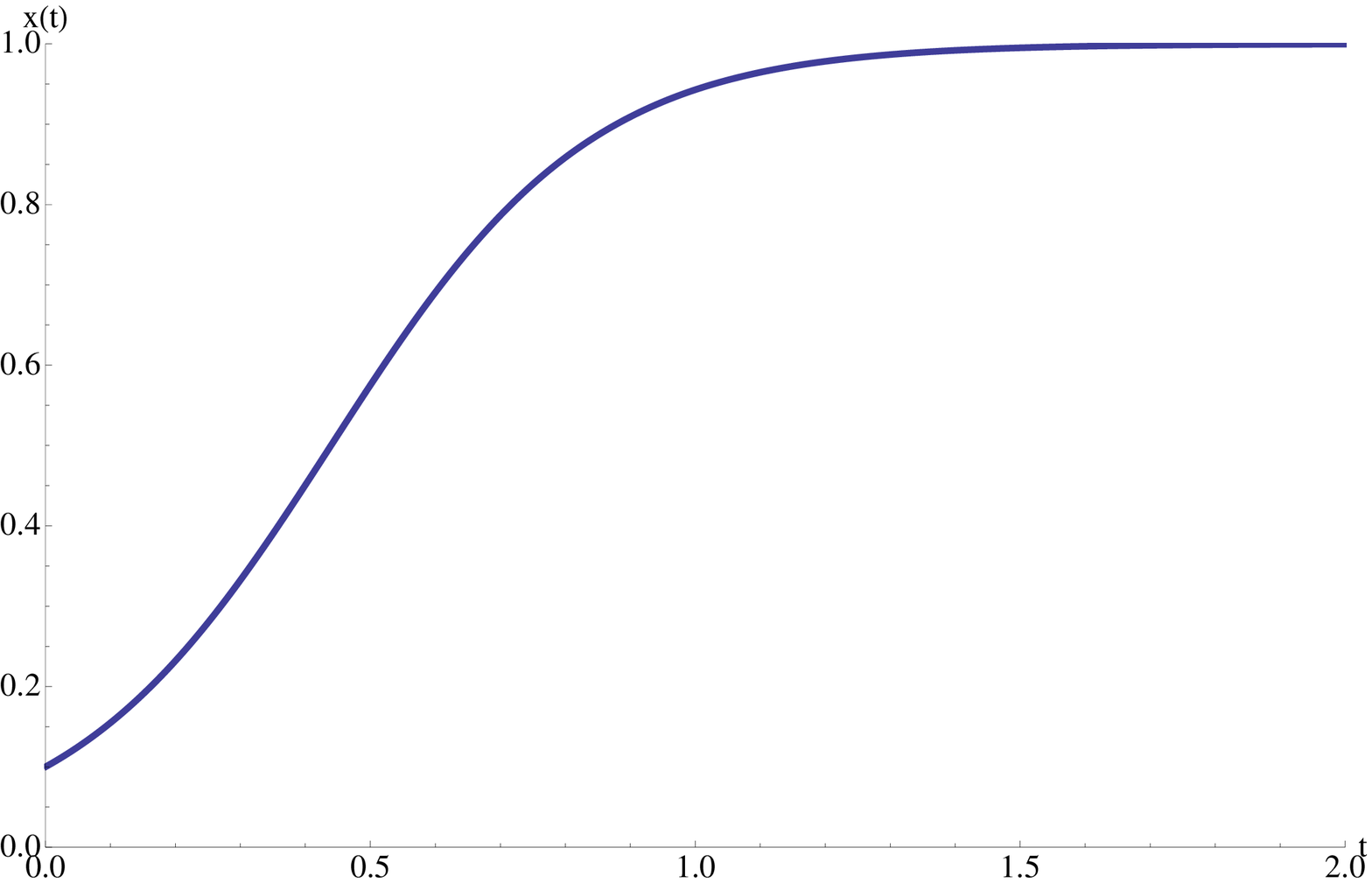}\label{fig: 1dode}}
\caption{Comparison of the CD methodology and ode45 for   the logistic equation.}\label{fig: mem}
\end{figure}

For this one-dimensional example, the initial value problem is
$$
Y'(t)= r Y(1-Y),\;\; Y(0) = Y_0
$$
with   the parameter $r =5$ and initial value $Y_0=0.1$. As we have already seen in the introduction, such problem is clearly non-convex with more than a local minimum even in its one dimensional version.

It is easy to notice that the matrices $\bG(\sigma_k)$ for this problem are scalars   equal to
$(\sigma_k + \sigma_{k+1}) r$ for $k=1,\dots,n$. The canonical dual variable $\bS = \{\sigma_k\} \in \real^n$ because of $d=1$.
Therefore $\inte \; \calS_a^+ \neq \emptyset$ as long as
  $  \sigma_k + \sigma_{k+1} \ge 0  \;\; \forall   k = 1, \dots, n$ and there exists one $\sigma_k  > 0$.
This nonlinear system clearly shows no chaotic behavior and the CD methodology is able to numerically find a global optimal solution
in $\calS_a^+$.

The trajectories by Runge-Kutta (ode45) and the CD methodology are reported in Figures \ref{fig: 1dcan} and \ref{fig: 1dode}.
The time interval in which the function has been analyzed is $[0,2]$ with step size $\delta = 1/500$, i.e. $n=1000$.
The target function (\ref{eq:primal}) for the two trajectories takes the following two values:
$$
P(\YY_{45})= 1.6141e-10, \;\; \quad P(\YY_{CD})=1.1219e-10.
$$
From these results it is clear that the two methods are able to give an accurate trajectory of the integral curve $Y(t)$ as the value of the objective function is close to zero. Because of the absence of chaotic behavior, no perturbation has been added to the problem.

\subsection{Forced Memristive Circuit  \cite{itoh11}}
We now turn our attention on a forced memristive circuit  governed by the following nonlinear system  \cite{chua76}:
\begin{equation}\label{eq: mem}
\begin{array}{rl}
\displaystyle \frac{dx}{dt}=&\displaystyle a (i^2-1),\vspace{0.5em}\\
\displaystyle\frac{di}{dt}=&\displaystyle-(x+\mu)i+\beta \cos (\omega t),
\end{array}
\end{equation}
where $a, \beta, \mu, \omega$ are given constants.
The input term  is the applied voltage source $v_s(t)=\beta \cos( \omega t)$.
Like  for the  previous example, we analyze the matrices $\bG(\bsigma_k)$ for $k=1,\dots,n$
 that are  2-dimensional   matrices:
$$
\left(
\begin{array}{cc}
0& -\frac{1}{2} (\sigma^2_{k}+\sigma^{2}_{k+1})\\
-\frac{1}{2} (\sigma^{2}_{k}+\sigma^2_{k+1}) & a(\sigma^{1}_{k}+\sigma^1_{k+1})
\end{array}
\right)  .
$$
The two  eigenvalues of this  matrix are
\eb
\begin{array}{rl}
\lambda_{1, 2}^k=&\frac{1}{2}a(\sigma^1_{k}+\sigma^1_{k+1}) \vspace{0.5em}\\
&\pm \frac{1}{2}\sqrt{a^2(\sigma^1_{k}+\sigma^1_{k+1})^2+(\sigma^2_{k}+\sigma^2_{k+1})^2 }.
\end{array}
\ee
Clearly,  the set $\calS^+_a$ is not empty only if   $\sigma^2_{  k} + \sigma^2_{  k+1} =  0$ for $k=1,\dots,n$. In this case,
$\inte \calS^+_a \neq \emptyset$ and the  problem  $(\calP^d)$ could have   solutions  on the boundary of $\calS^+_a$.
The  conventional iterative methods could produce pseudo-chaos,
for example, in the case of  $x_0=0.1,i_0=0.1,a=1,\mu=0,\beta=0.7,\omega=1$.
But by using the perturbed canonical dual methodology, we are still able to find global optimal trajectory even for this
chaotic system.
 Table \ref{tab:mem} shows the global errors $\Pi(\bY)$ produced by R-K method (ode45) and perturbed CD method with different step sizes $\delta$,
 initial values $Y_0$, amplitudes $\beta$,
 and perturbation $\rho$.

\begin{table}[ht]
\begin{small}
\begin{center}
\begin{tabular}{|c| c| c| c|c| c|c|}
\hline
n & $Y_0$&$\beta$&	$\Pi(\bY_{45})$		&	$\rho_0$	&	$\rho^*$ &	$\Pi(\bY_{CD})$	\\
\hline														
1000	& (0.1,0.1)&0.7&	6.4434	&	1	&	$1.3841* 10^{-2}$	&	2.8172	\\	\hline
10000& (0.1,0.1)&0.7	&	0.6237	&	$10^{-2}$	&	$3.125*10^{-4}$	&	0.334\\	\hline
100000& (0.1,0.1)&0.7	&	0.0624	&	$10^{-3}$	&	$7.29* 10^{-7}$	&	0.0275\\	
\hline
100000& (0.2,0.2)&1&0.0751&$10^{-4}$&$8.1*10^{-7}$&0.0223\\
\hline
\end{tabular}
\end{center}
\caption{Global errors $P(\bY)$ produced by R-K (ode45)  and CD methods  for the forced memristive circuit. }
\vspace{0.5em}
\label{tab:mem}
\vspace{-.5cm}
\end{small}
\end {table}


\begin{figure*}
\subfigure[CD  solution   $x(t)$.]{\includegraphics[scale=.13]{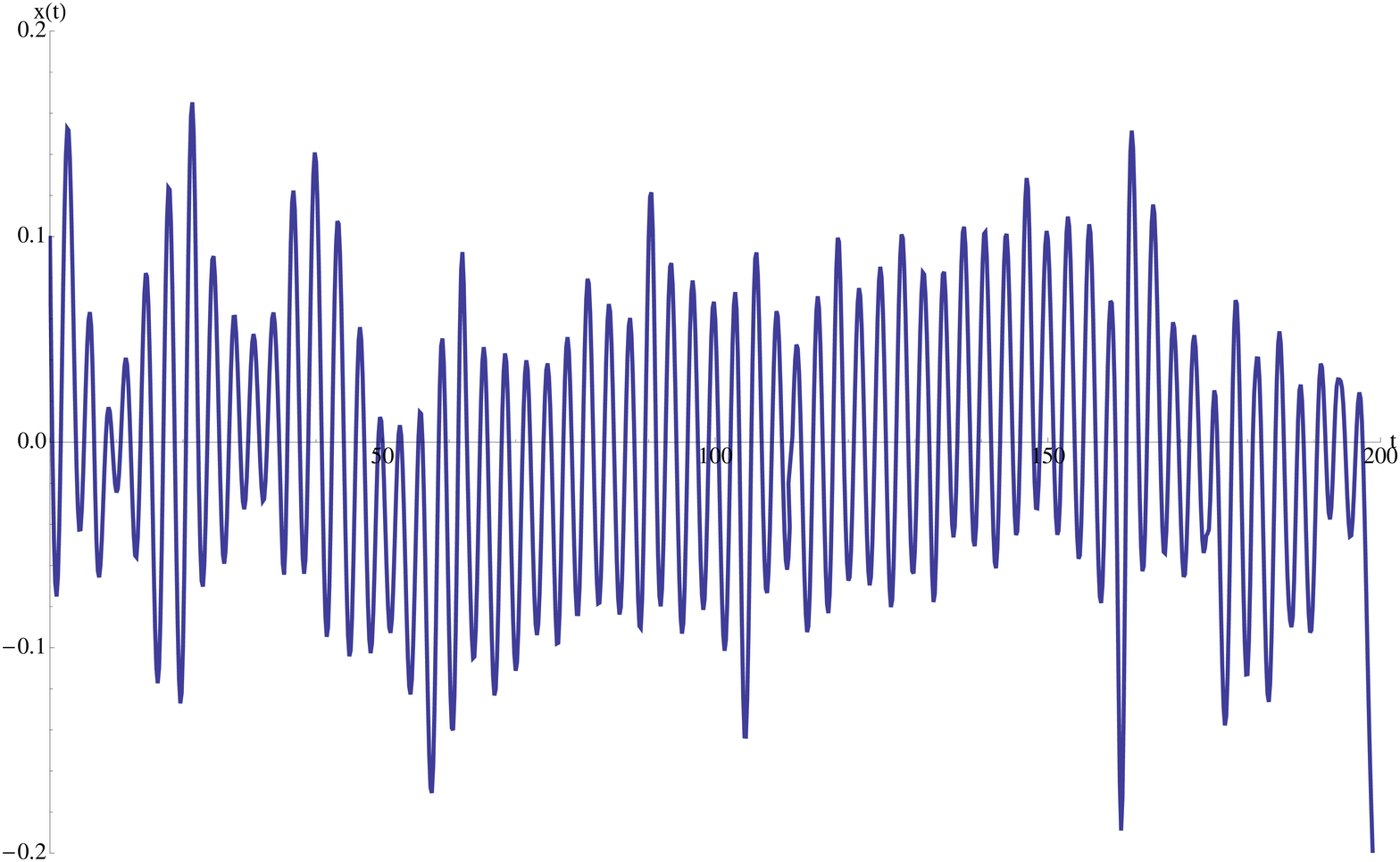}}
\subfigure[CD solution  $i(t)$.]{\includegraphics[scale=.13]{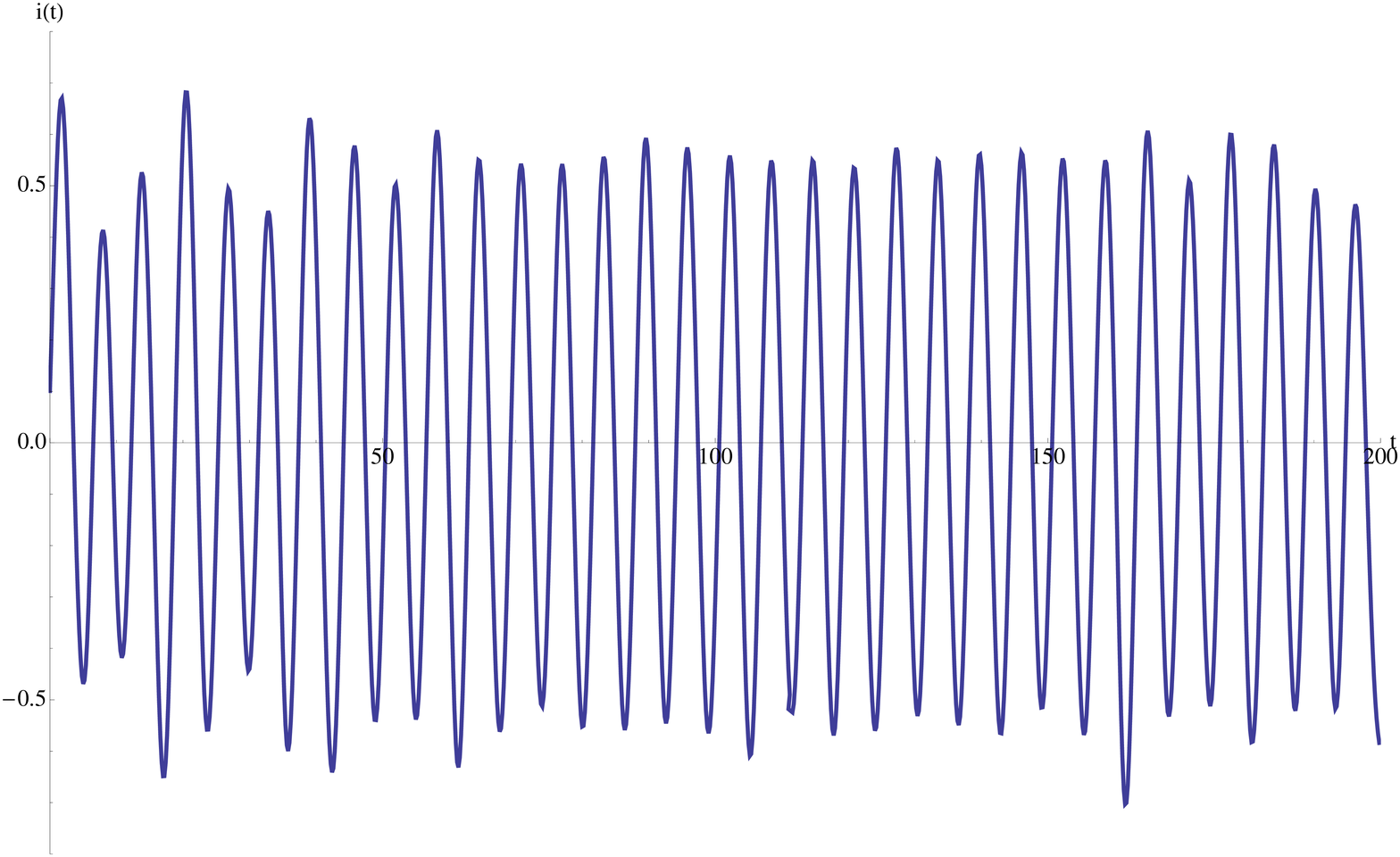}\label{memcany}}
\subfigure[Trajectory  by CD methodology.]{\includegraphics[scale=.13]{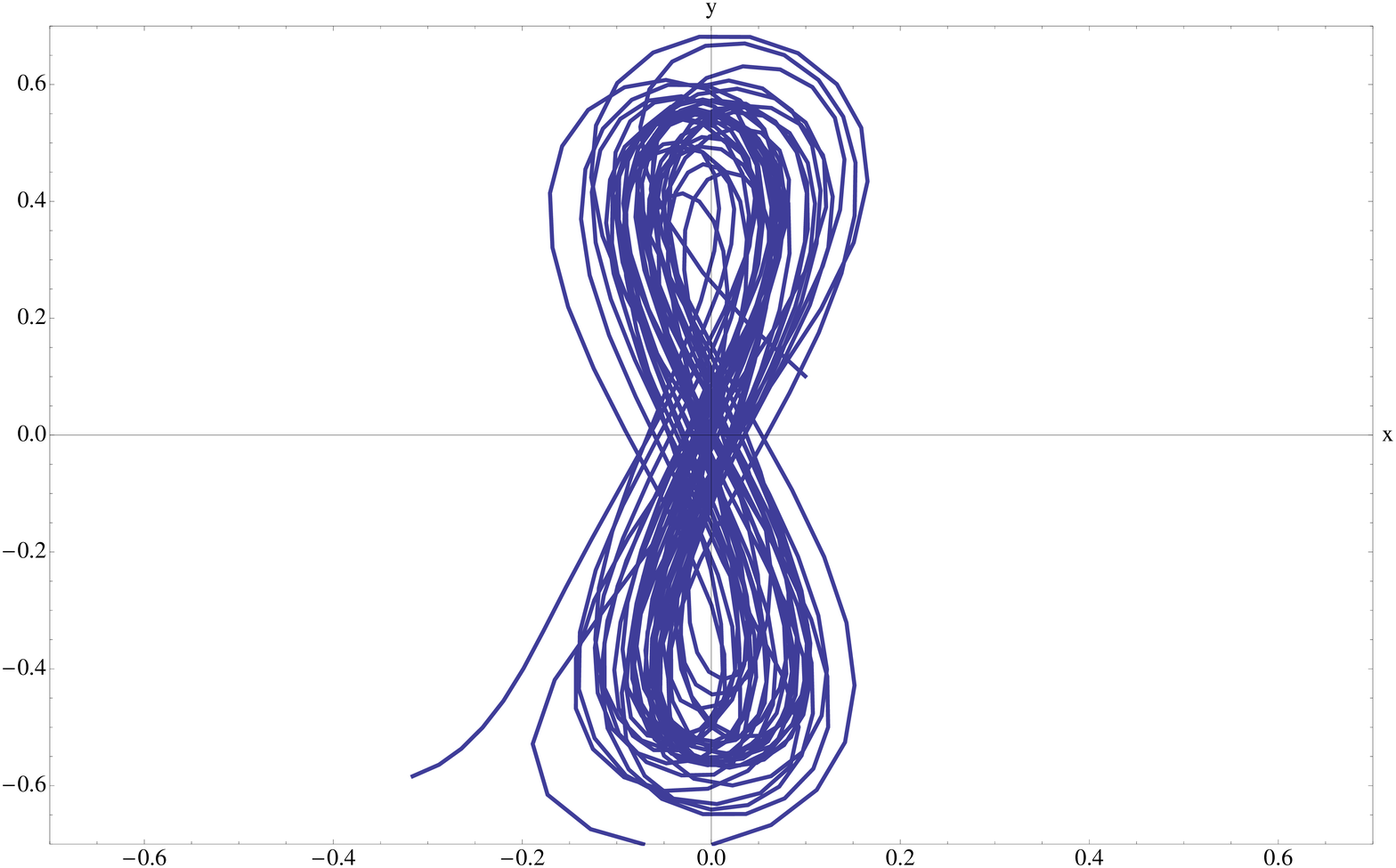}\label{memcanxy}}

\subfigure[R-K solution   $x(t)$ by ode45.]{\includegraphics[scale=.13]{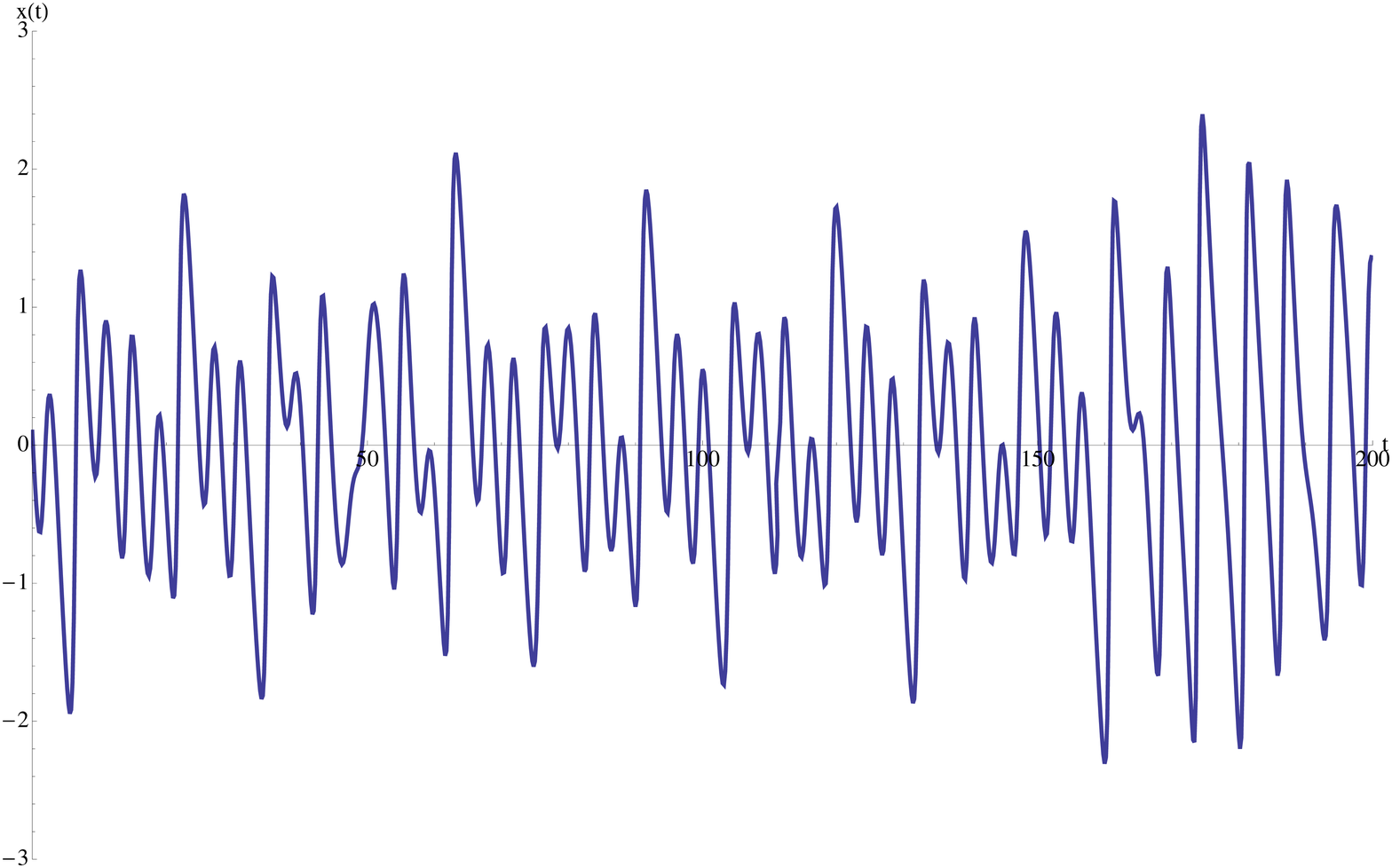}}
\subfigure[R-K solution  $i(t)$  by ode45.]{ \includegraphics[scale=.13]{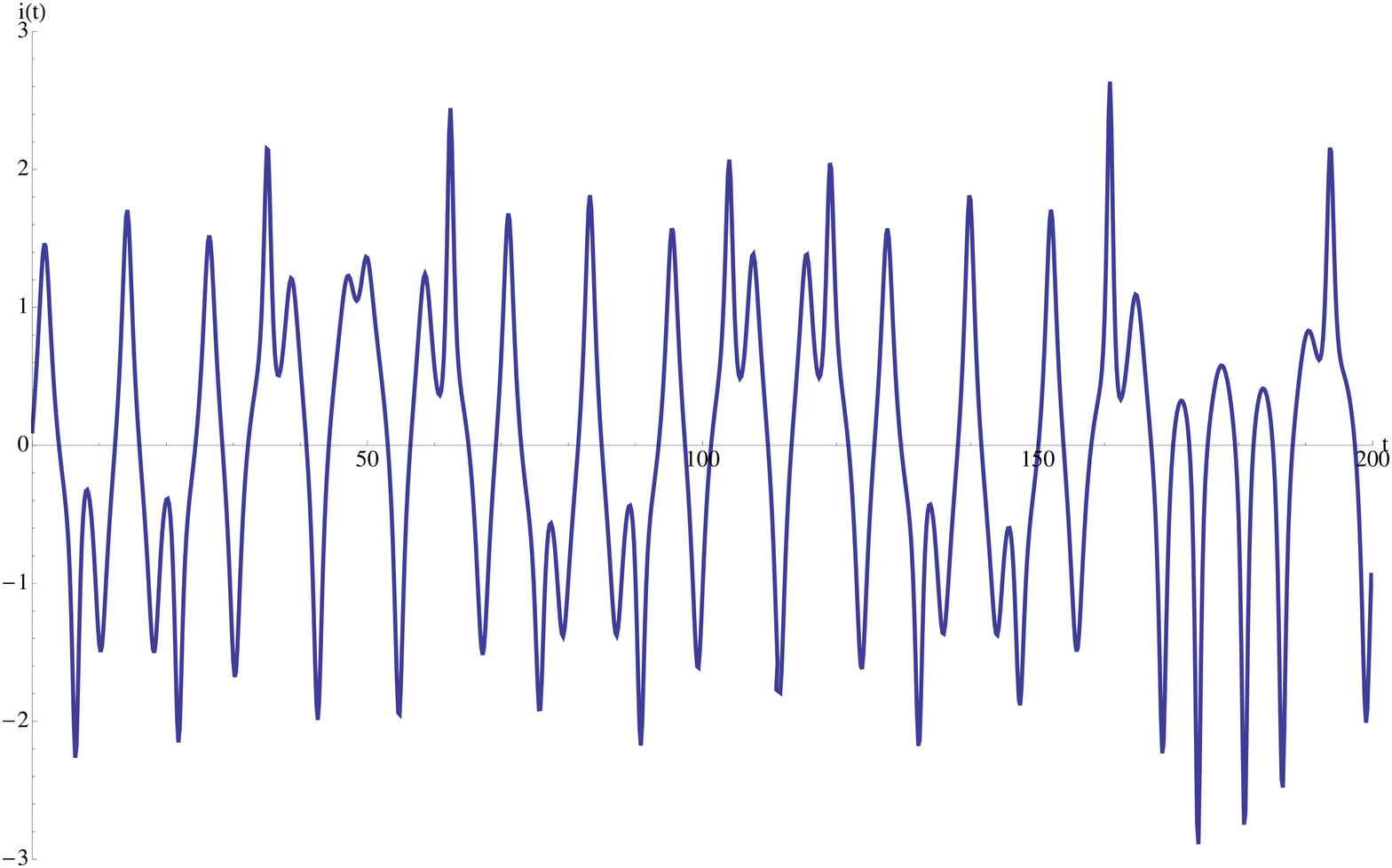}\label{memodey}}
\subfigure[Trajectory   by  ode45.]{\includegraphics[scale=.13]{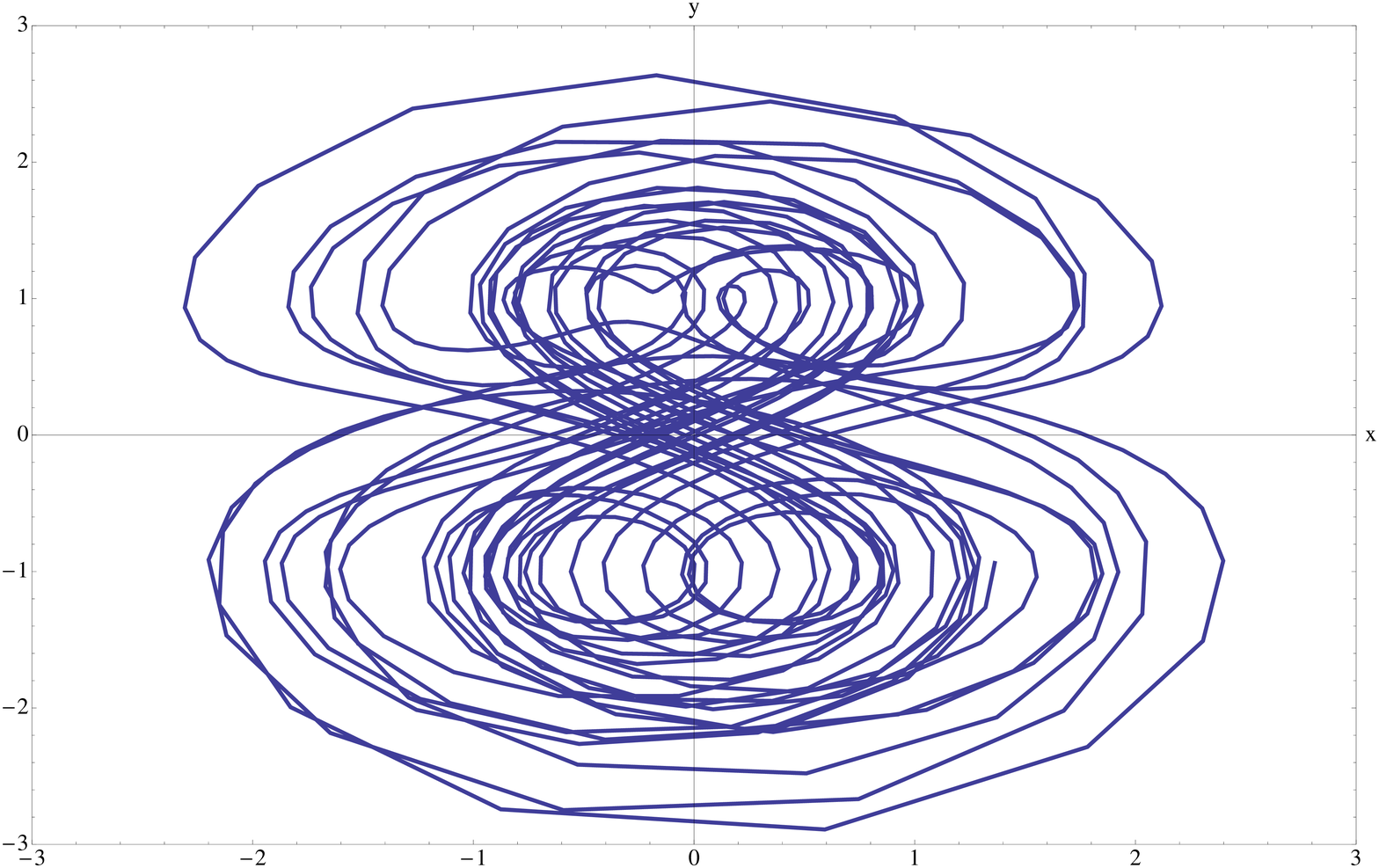}\label{memodexy}}
\caption{Comparisons of the CD methodology and R-K method (ode45) for   Eqn (\ref{eq: mem}) with $n = 100,000$, $\beta = 0.7$, and $Y_0 =( 0.1, 0.1 )$}\label{fig: mem}
\end{figure*}

 Figure \ref{fig: mem}  shows the chaotic trajectories produced by conventional R-K iteration method with
 oscillations between $(-3, 3)$, see Figure \ref{fig: mem} (d).
While  the global optimal solutions produced by the perturbed canonical dual method are   quite more regular and stable with
 oscillations between $(-0.7, 0.7)$, see Figure \ref{fig: mem} (c).

The results  with changed initial point and amplitude are reported in Figure \ref{fig: smem}.
Compared with  Figure \ref{fig: mem}  we can see that the Runge-Kutta method still   produces the  pseudo-chaotic ``solutions",
while the canonical duality method produces more stable and realistic global optimal  solutions.




\begin{figure*}
\subfigure[CD solution   $x(t)$.]{\includegraphics[scale=.13]{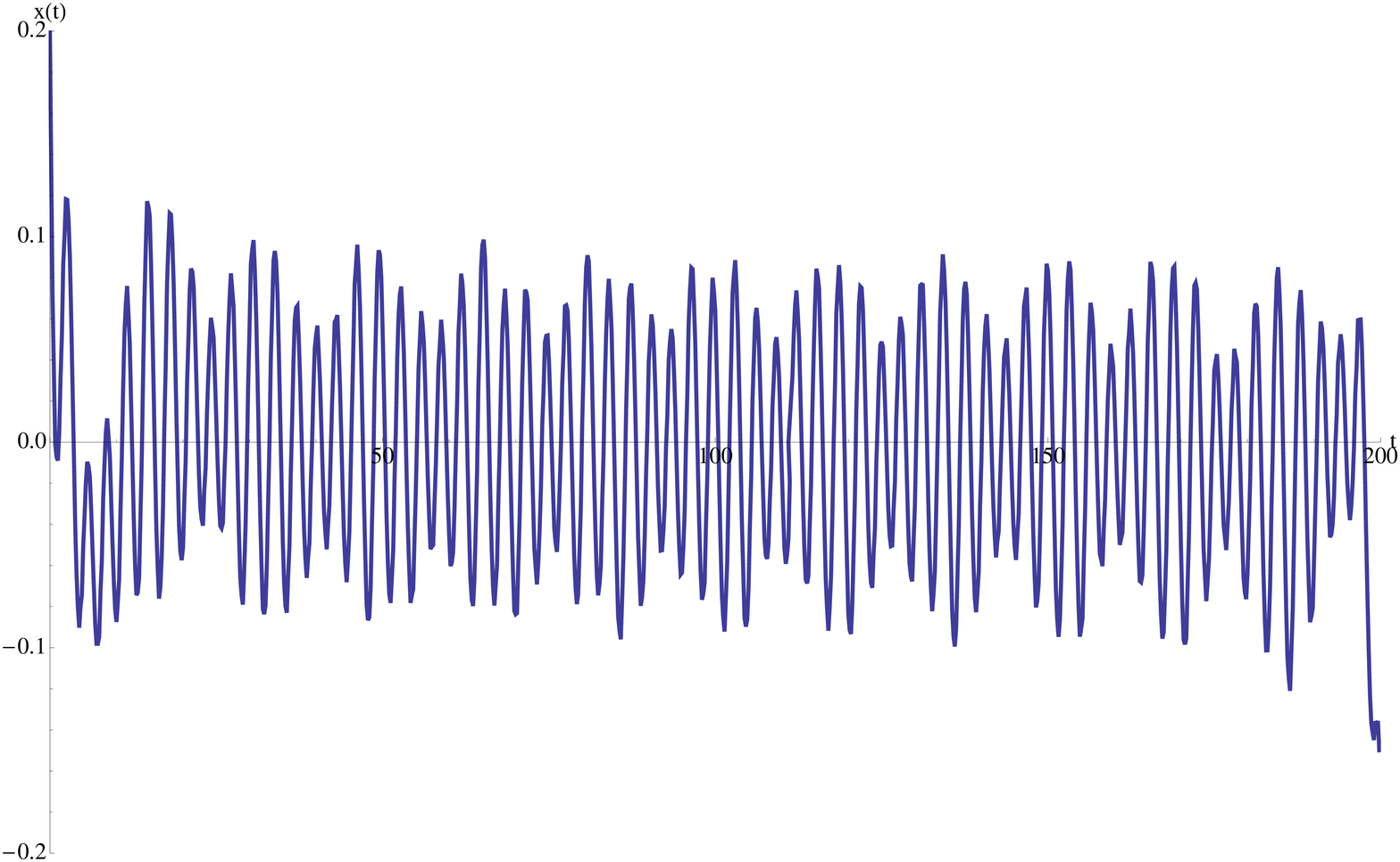}}
\subfigure[CD  solution  $i(t)$.]{\includegraphics[scale=.13]{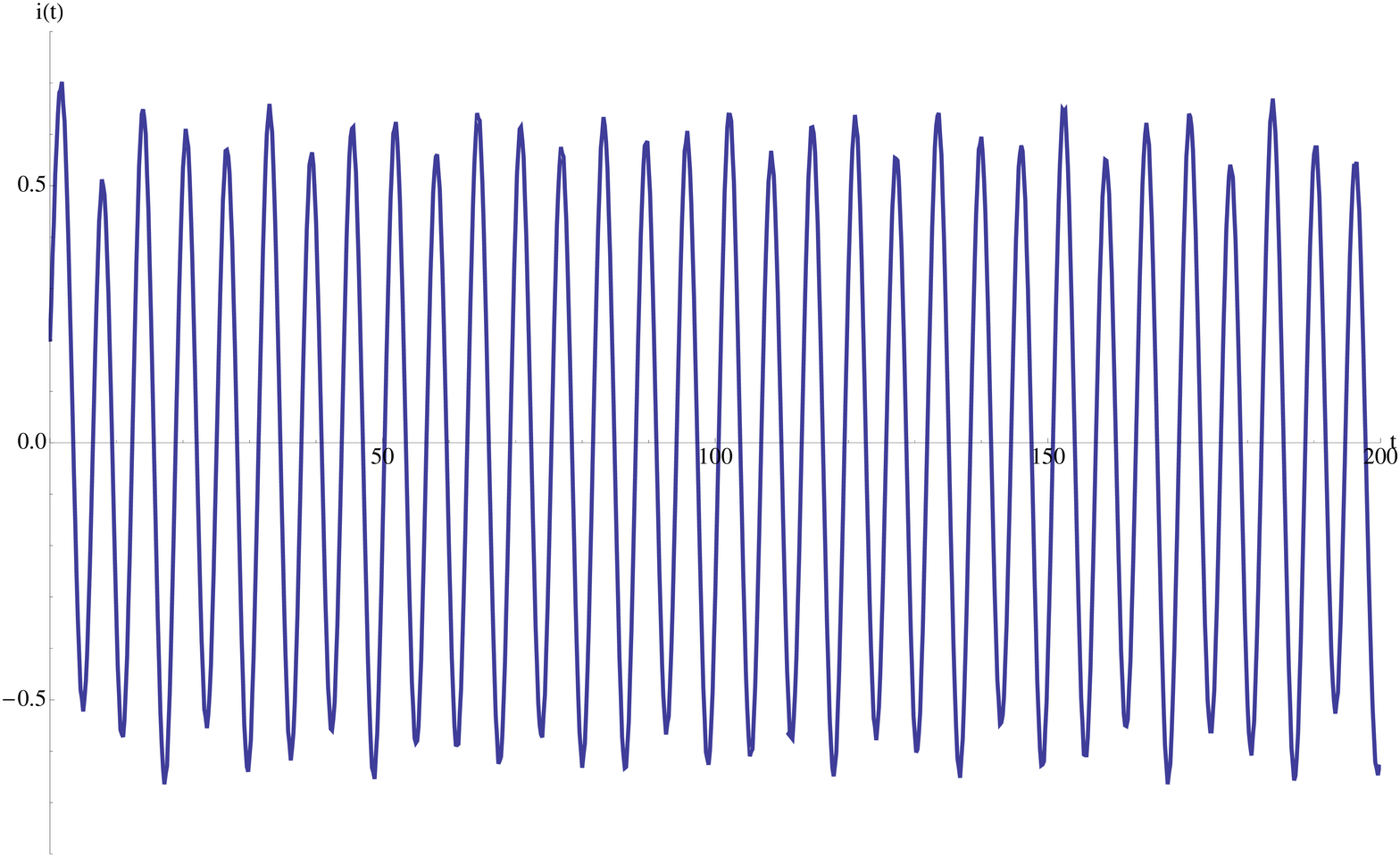}\label{smemcany}}
\subfigure[CD trajectory.]{\includegraphics[scale=.13]{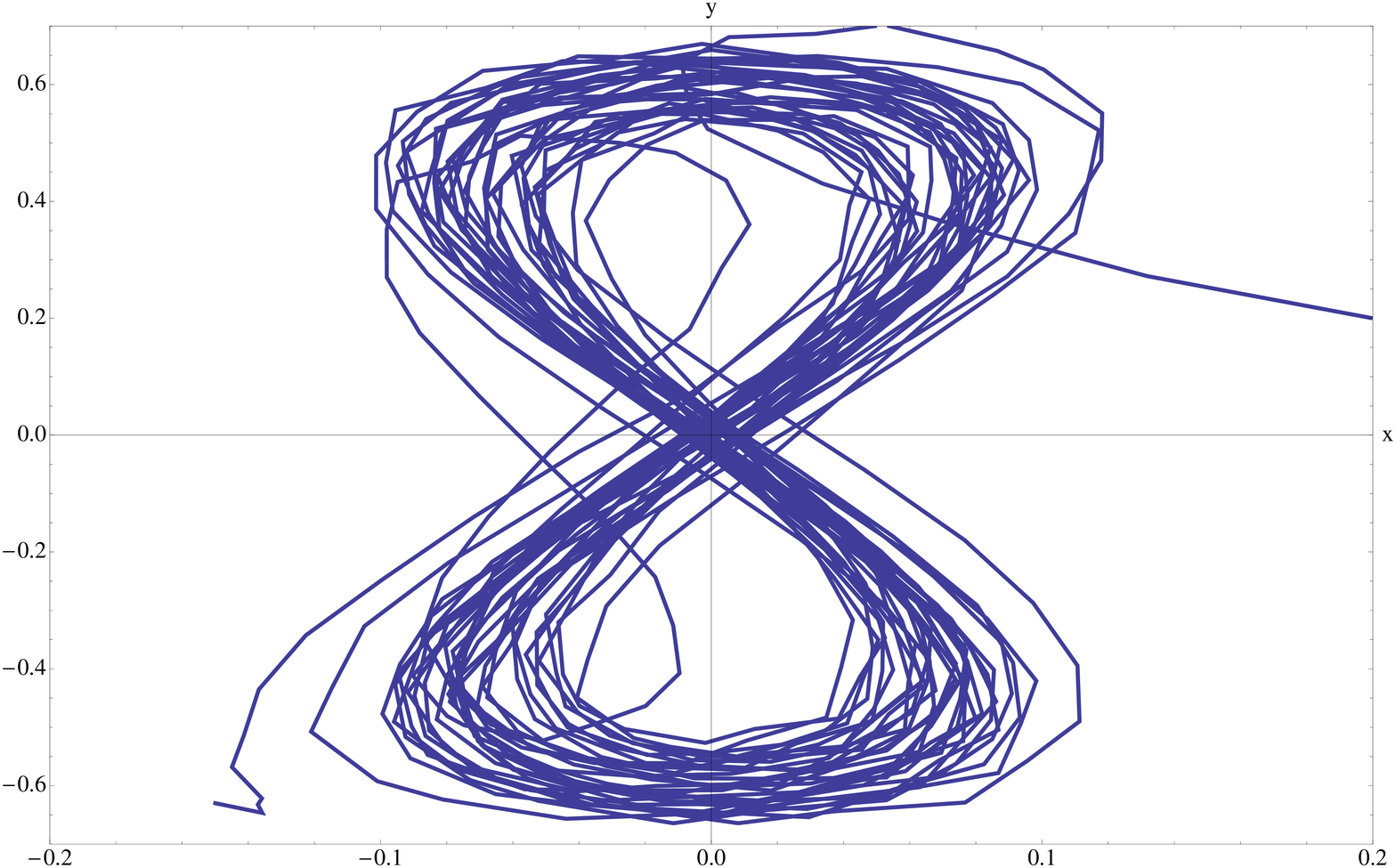}\label{smemcanxy}}

\subfigure[R-K solution  $x(t)$  by ode45.]{\includegraphics[scale=.13]{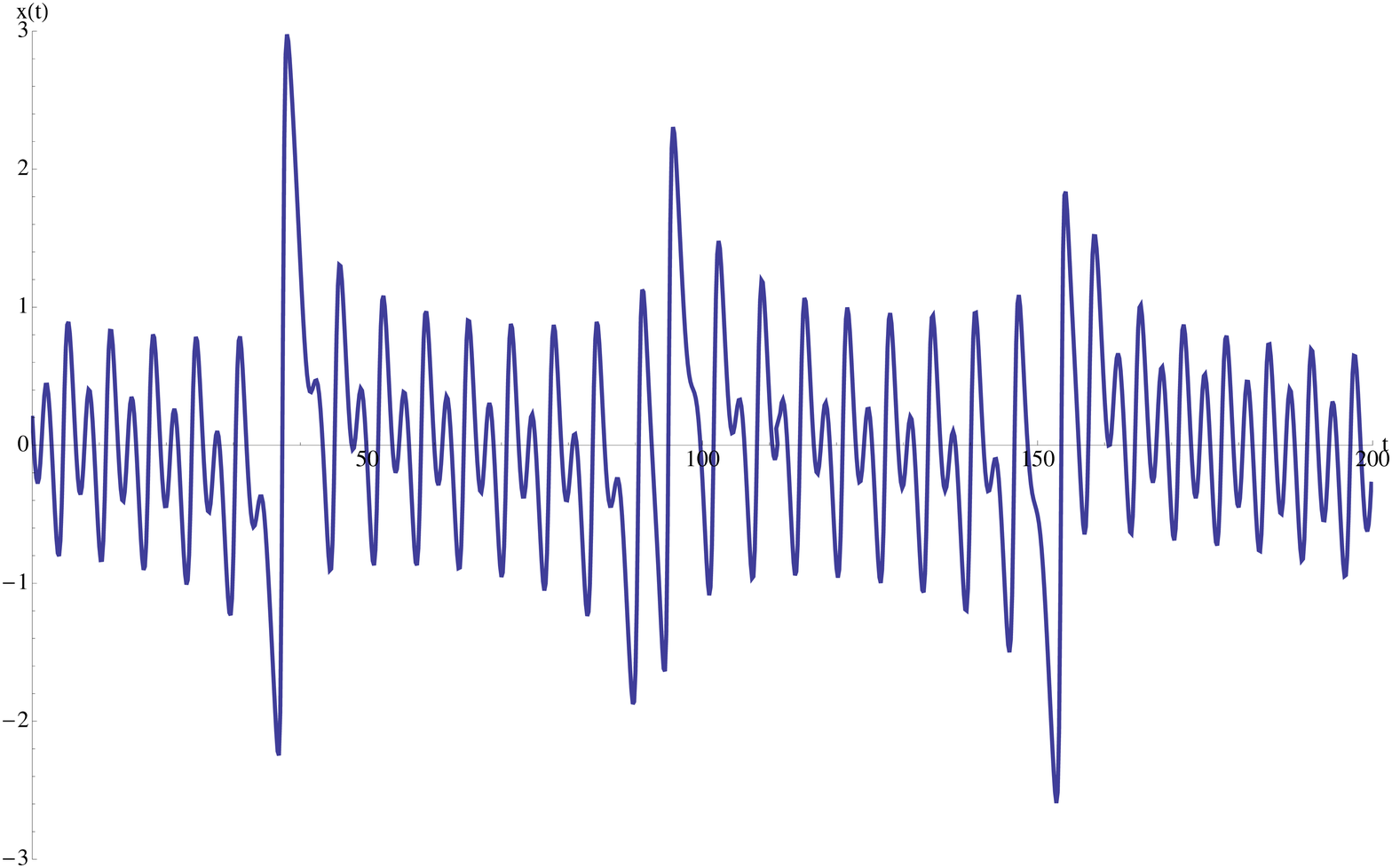}}
\subfigure[R-K  solution of $i(t)$  by  ode45.]{ \includegraphics[scale=.13]{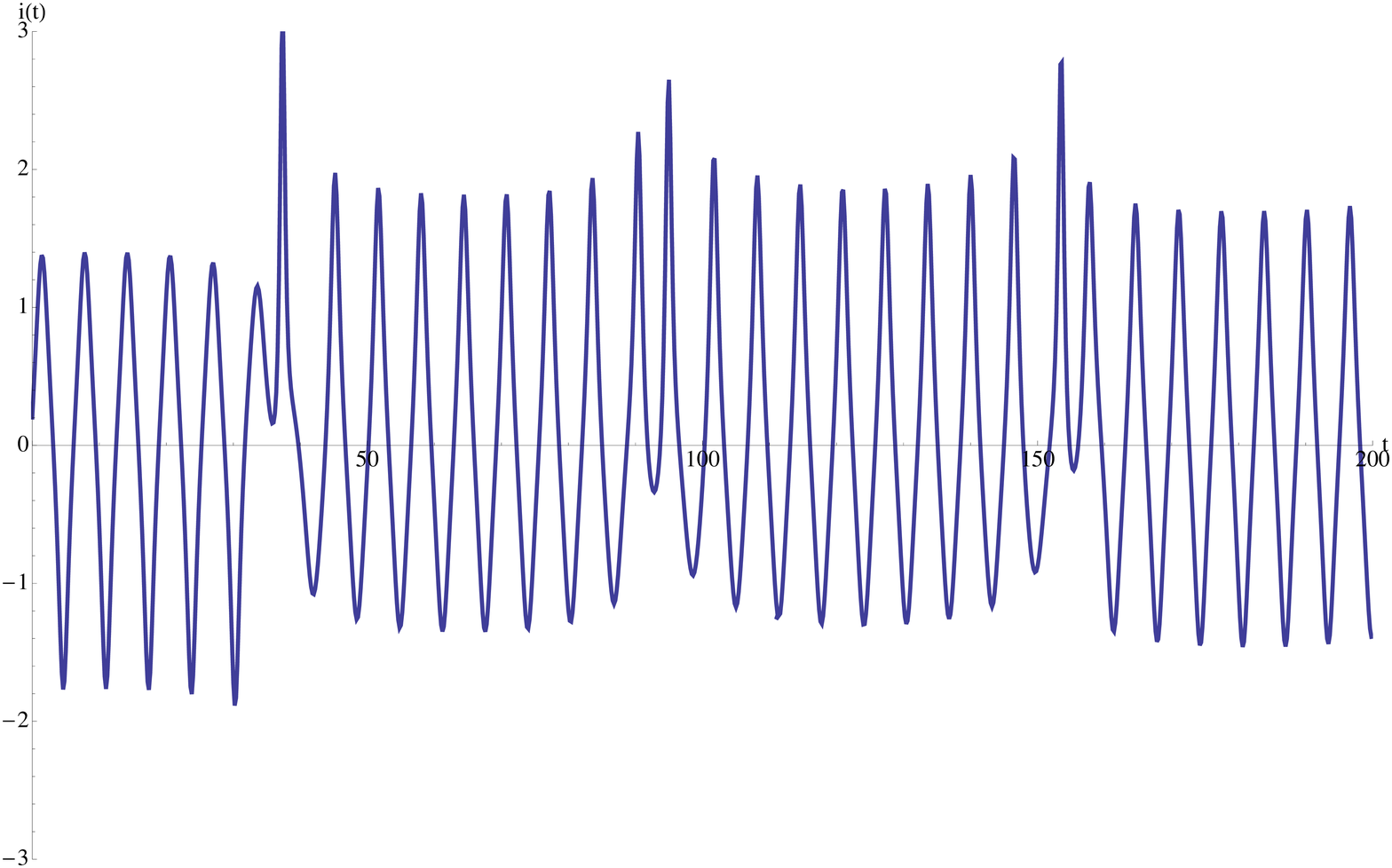}\label{smemodey}}
\subfigure[R-K trajectory  by  ode45.]{\includegraphics[scale=.13]{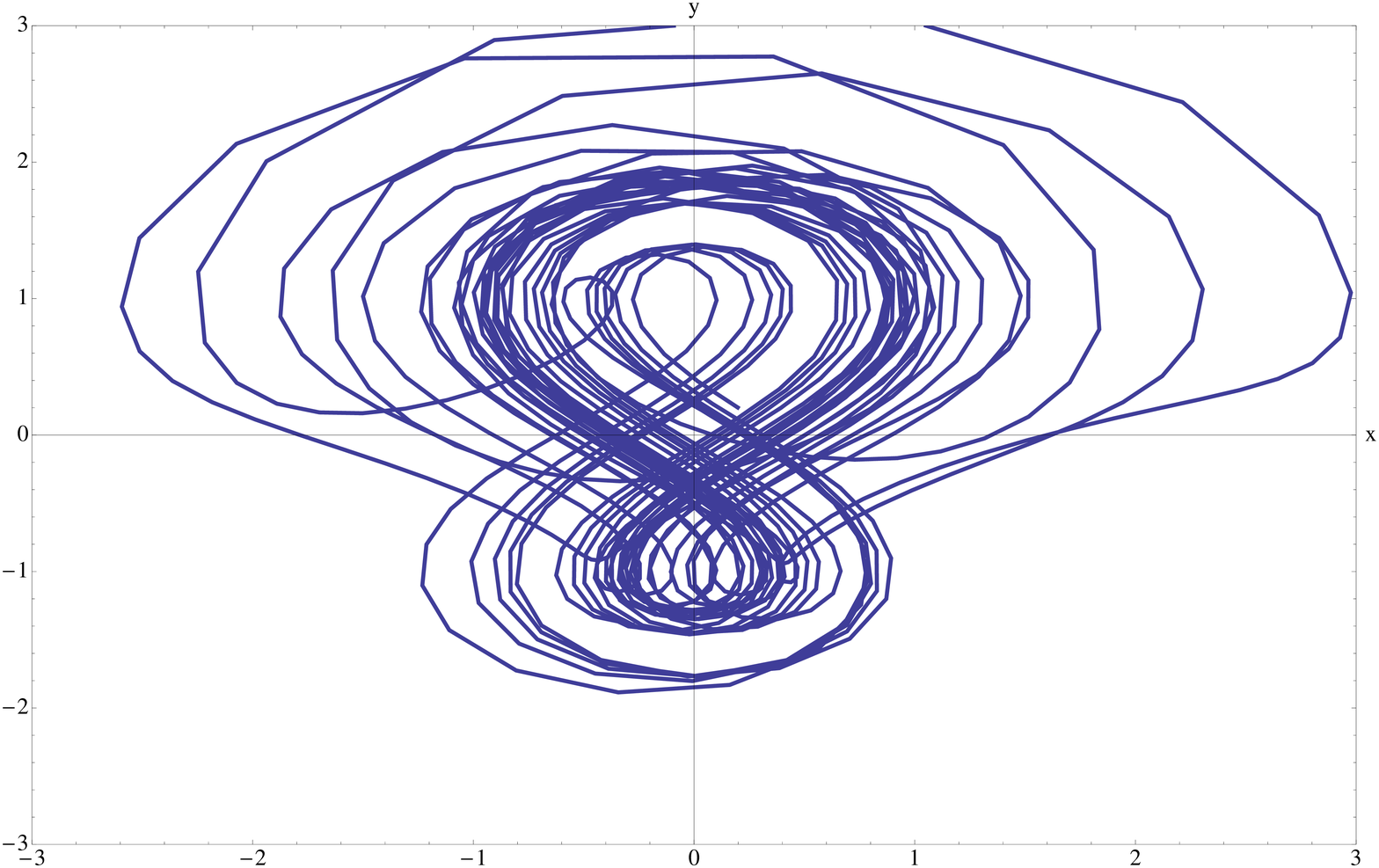}\label{smemodexy}}
\caption{Comparisons of the CD methodology and R-K method for   Eqn (\ref{eq: mem})  with $n=100,000$, $Y_0=(0.2,0.2)$,
 $\beta=1$}\label{fig: smem}
\end{figure*}

\begin{figure}
\subfigure[R-K solutions  $x(t)$  by   ode45 (blue) and ode23 (red).]{\includegraphics[scale=.095]{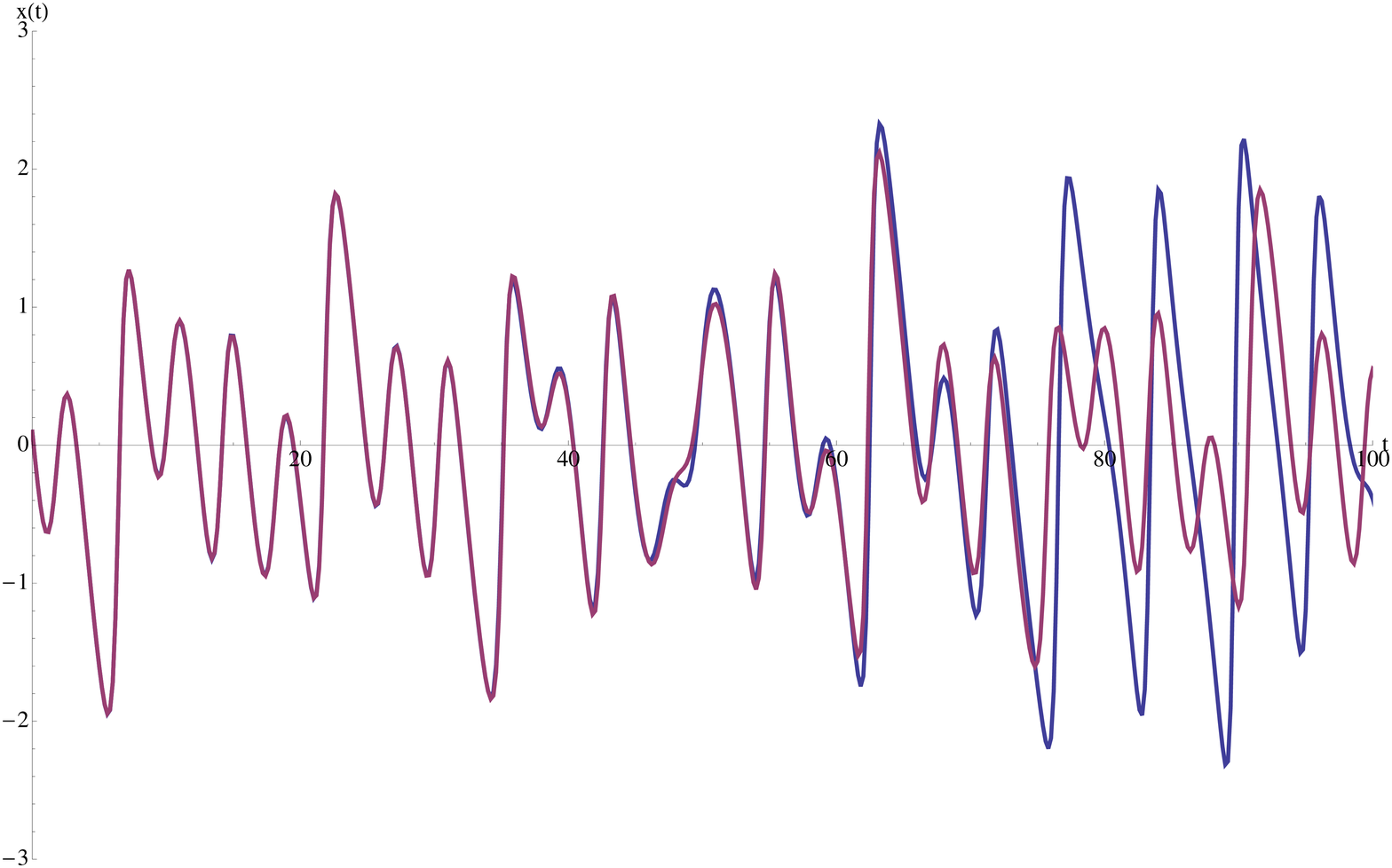}} \hspace{.3cm}
\subfigure[R-K solutions  $i(t)$ by  ode45 (blue) and ode23 (red).]{\includegraphics[scale=.095]{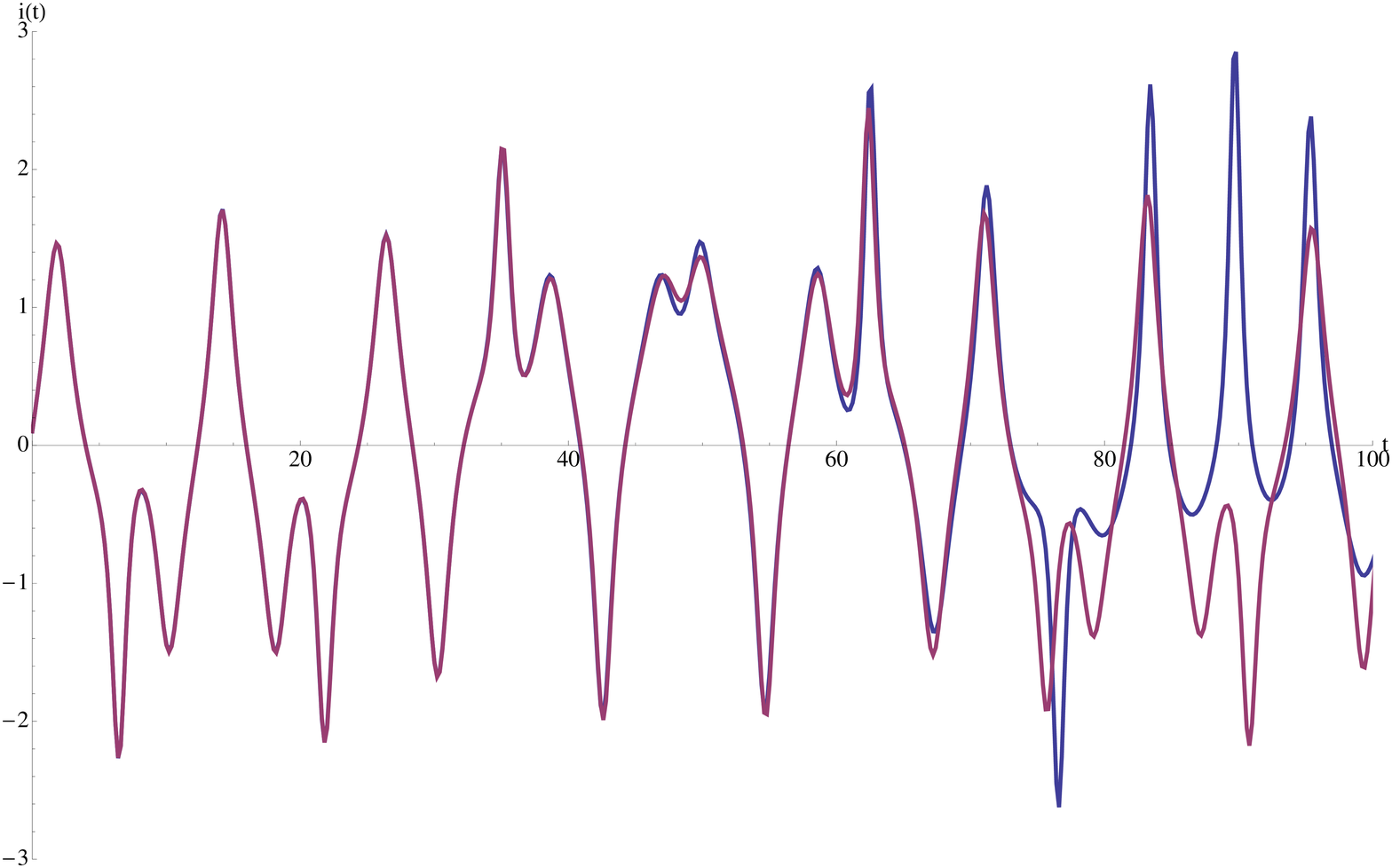}}
\caption{Comparisons of R-K solutions for Eqn  (\ref{eq: mem})
produced by  ode45 and ode23 with $n=100000$ for the memristive circuit. It is possible to see that the two trajectories are diverging  with the increasing number of iterations }\label{fig: smem1}
\end{figure}
To further understand the {\em pseudo-chaos}  due to  error-accumulations produced by  the conventional iteration methods,
  we present the chaotic numerical results to the  Eqn  (\ref{eq: mem}) produced by  an another  routine {\em ode23 } in
  Matlab, with is also based on Runge-Kutta method.
  The total error of this solution is $\Pi(\bY_{23})= 0.0622$. From  Figure \ref{fig: smem1}  we can see that the two trajectories
   coincide for the first iterations and then, as the numerical error  accumulates, they diverge.
This fact shows that the conventional
   iterative methods can't produce  reliable numerical solutions for deterministically stable
    nonlinear systems due to the intrinsic numerical error accumulations.

\newpage

\subsection{Lorenz Equations}
The Lorenz system  has been  studied extensively in literature and is  well known for its highly chaotic behavior with
butterfly attractor \cite{lo12}.
This section presents a canonical dual approach to show the reason for chaos and to verify Conjecture \ref{con: c=np}.

The initial value problem for Lorenz system  with standard parameters is
$$
\begin{array}{rl}
\displaystyle\frac{dx}{dt}&\displaystyle=10(y-x) , \;\; x(0)  = 10 , \vspace{0.5em}\\
\displaystyle\frac{dy}{dt}&\displaystyle=-xz+28x-y  , \;\; y(0) = 12, \vspace{0.5em}\\
\displaystyle\frac{dz}{dt}&=\displaystyle xy-\frac{8}{3}z , \;\; z(0) = 14 .
\end{array}
$$
In order to understand how $\calS_a^+$ is characterized, we start the analysis by looking at the matrices $\bG(\bsigma_k)$ for $k=1,\dots, n$:
$$
\left(
\begin{array}{ccc}
0&\frac{1}{2}(\sigma^3_{k}+\sigma^3_{k+1}) &-\frac{1}{2} (\sigma^2_{k}+\sigma^2_{k+1})\\
\frac{1}{2}(\sigma^3_{k}+\sigma^3_{k+1})&0&0\\
-\frac{1}{2}(\sigma^2_{k}+\sigma^2_{k+1}) & 0&0
\end{array}
\right) .
$$
The  three eigenvalues can be easily calculated as:
\begin{equation}
\begin{array}{rl}\label{eq: loeig}
\displaystyle\lambda_1^k&\displaystyle=0\vspace{0.5em},\\
\displaystyle\lambda_2^k&\displaystyle=\frac{1}{2}\sqrt{(\sigma^2_{k}+\sigma^2_{k+1})^2+(\sigma^2_{k}+\sigma^2_{k+1})^2}\vspace{0.5em},\\
\displaystyle\lambda_3^k&\displaystyle=-\frac{1}{2}\sqrt{(\sigma^2_{k}+\sigma^2_{k+1})^2+(\sigma^2_{k}+\sigma^2_{k+1})^2}.\\
\end{array}
\end{equation}
Equations (\ref{eq: loeig}) show that for any  $\bS = \{\bsigma_k \} \in \calE^*_a$, the global  matrix $\{\bG(\bsigma_k)\} $  has total $n$ zero eigenvalues, $n$ positive eigenvalues  and $n$ negative eigenvalues.
Therefore   $ \bG(\bsigma_k)  $ is  indefinite for any   dual
 variable  $\bS = \{ \bsigma_k\} \in \calE^*_a $, thus both $\inte \; \calS_a^+$ and $\calS^-_a $ are  empty for  the Lorenz system.
By Conjecture \ref{np-hard} we know that the problem $(\calP)$ is NP-hard.
This explains the highly chaotic behavior of the Lorenz system and verifies Conjecture \ref{con: c=np}.
\begin{figure}
\subfigure[CD solution  $x(t)$  using initial iteration produced by  {ode45}.]{\includegraphics[scale=.1]{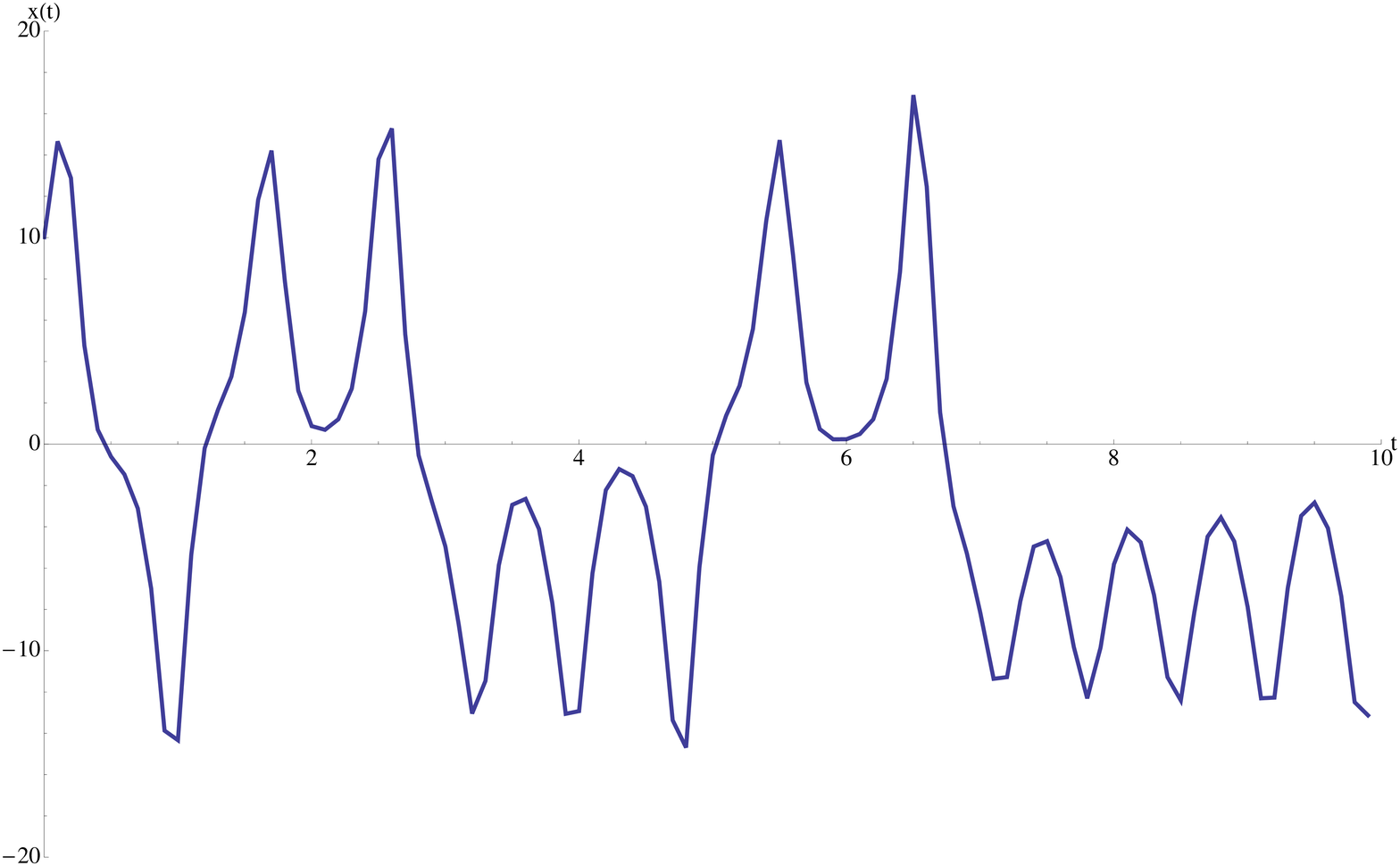}\label{4locanx}}
\subfigure[R-K solution  $x(t)$   by ode45.]{\includegraphics[scale=.1]{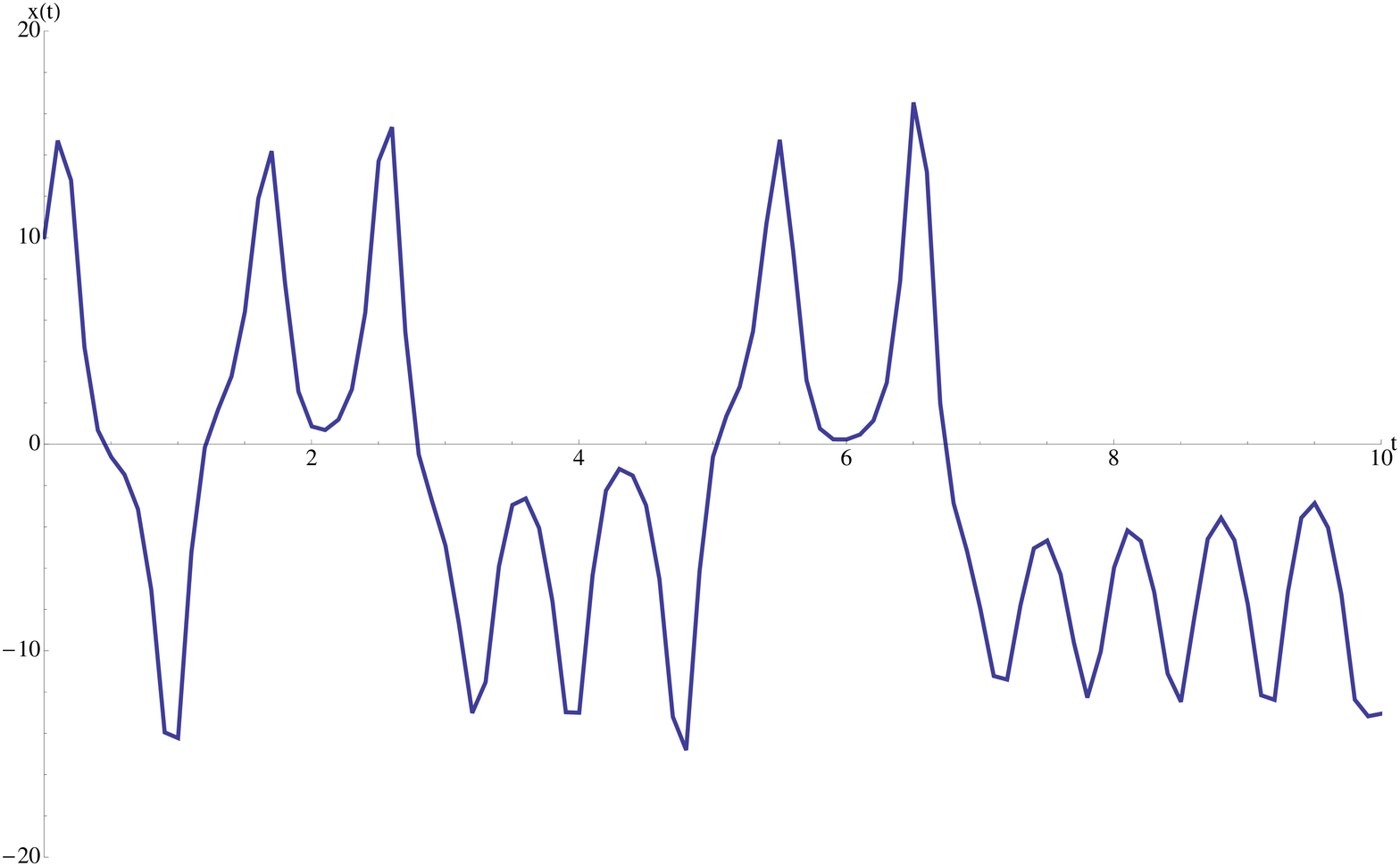}\label{4lodex}}
\caption{Comparison between the CD methodology and ode45 for the Lorenz system}\label{fig: lo}
\end{figure}

By the fact that $\inte \; \calS^+_a = \emptyset$ for the Lorenz system,
 we have to solve the nonconvex canonical dual problem $(\calP^d_s)$ defined by (\ref{cd-s}),
  which is NP-hard and numerical solutions depend sensitively on the initial iteration.
 We have  performed the  test  for
 time interval  $[0,10]$ with  step size $\delta = 1/1000$, i.e. $n=10,000$. First we run ode45 and then use the numerical result by Runge-Kutta (R-K solution) as a starting point for the canonical dual iteration. The numerical solutions for $x(t)$ are reported in Figures \ref{4locanx} and \ref{4lodex}. It is possible to see that the two trajectories are quite similar, but
 the main difference is the total errors:
$$
\Pi (\bY_{45})=10^{-3},\;\;\; \quad  \Pi (\bY_{CD})=3.326*10^{-07},
$$
i.e., the total error produced by CD method is much smaller than the one produced by R-K method.
 By the fact that   $\calS_a^+$ has no-interior, we can't theoretically claim that $\bY_{CD}$ is a
global minimum solution to Problem (\ref{eq:primal}) for the Lorenz equations, however,  the fact that $\Pi (\bY_{CD}) \simeq 0$  indicates
 that the canonical dual solution  $\bY_{CD}$  is indeed one of  global optimal solutions to
the nonlinear system $(\calP_0)$.
By Conjecture \ref{con: c=np} we can say  that the Lorenz system is chaotic and the CD method produces the well-known  butterfly trajectories as  shown in Figure \ref{butterfly}.
\begin{figure}
\subfigure[ Trajectory  produced by  ode45.]{\includegraphics[scale=.095]{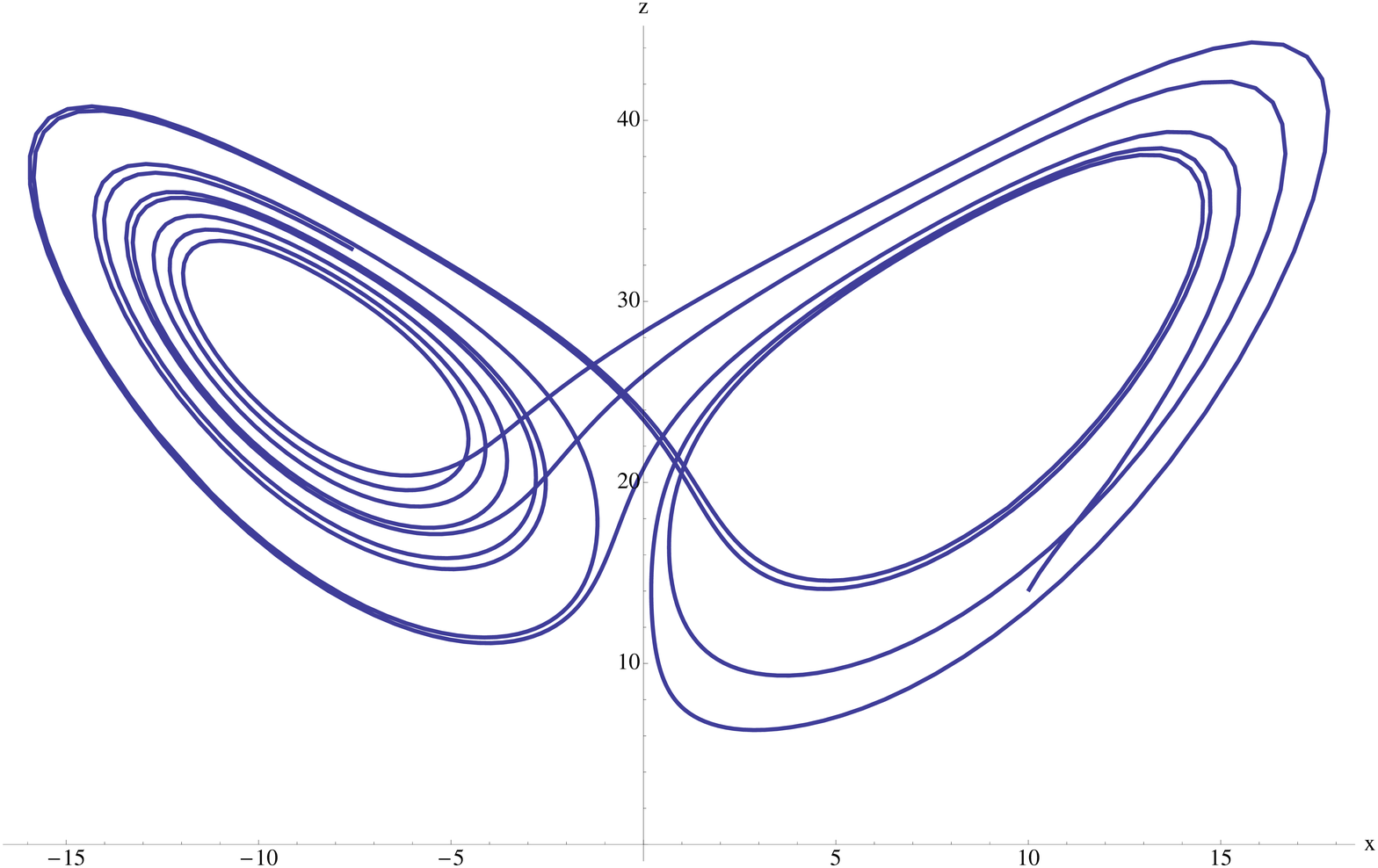}}
\subfigure[Trajectory produced by CD method.]{\includegraphics[scale=.095]{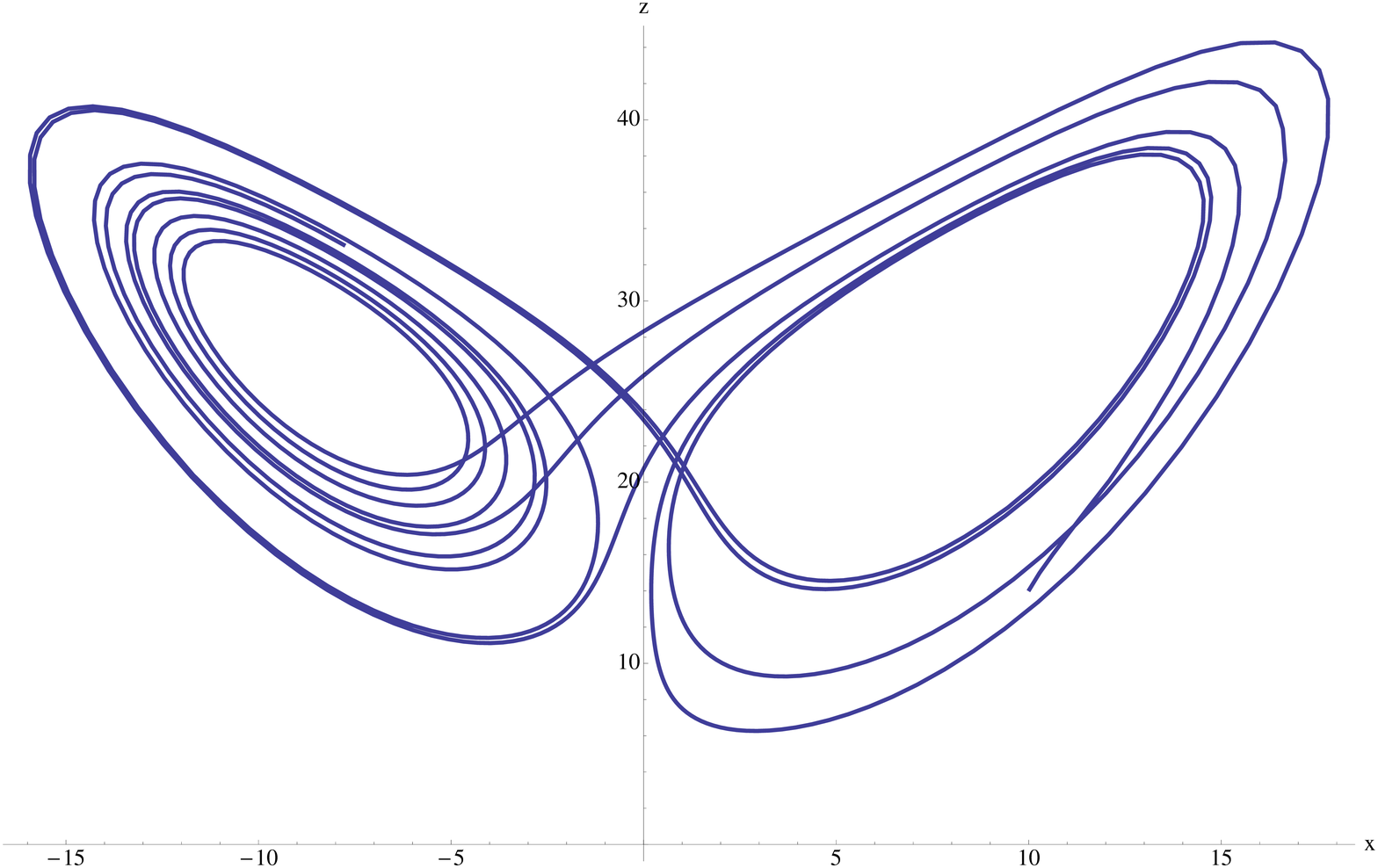}}
\caption{Lorenz attractors produced by  ode45 (left) and  the CD methodology (right)}\label{butterfly}.
\end{figure}

Finally in Figure \ref{locanx} we report the trajectories produced  by the perturbed CD methodology starting from the origin for the values of the primal variables. In detail:
\begin{itemize}
\item the Yellow curve represents the solutions for $x(t)$ with $\rho=10^{-3}$ and $\Pi_\rho(\bY^*)=4.9$;
\item The red curve represents the solutions for $x(t)$ with $\rho=10^{-5}$ and $\Pi_\rho(\bY^*)=2.5$;
\item The blue curve represents the solutions for $x(t)$ with $\rho=10^{-6}$ and $\Pi_\rho(\bY^*)=1.6$;
\end{itemize}
We can see that all of these trajectories are stable even after changing the value of $\rho$. Therefore we can state that if the algorithm by the CD methodology starts from a generic point in $\calS_\rho^+ $ with an high enough value of $\rho$ it  reaches a stable trajectory. On the other side, because of the characterization of the original $\calS_a^+$, the error is quite high, far from the theoretical optimal value that is zero. Nevertheless such trajectories correspond to local solutions of Problem (\ref{eq:primal}).
\begin{figure}
 {\includegraphics[scale=.2]{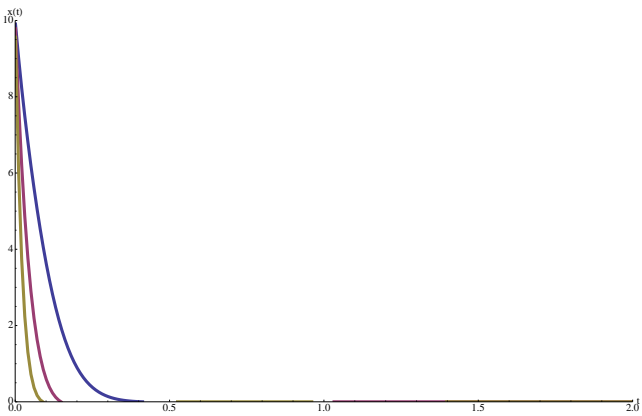}
\caption{Perturbed CD  solutions  $x(t)$   starting from the origin with $\rho= 10^{-6}$ (blue),  $\rho= 10^{-5}$ (red), and  $\rho =10^{-3}$ (yellow).}\label{locanx}}
\end{figure}

 \section{Conclusions}

 We have presented a nonconventional method for solving  nonlinear dynamical systems.
   Based on the  methods of finite differences and least squares, the nonlinear initial value problem is first reformulated as  a
  global minimization problem.
  By using the canonical duality theory, the nonconvex minimization problem is able to converted to a
     concave maximization dual problem over a convex cone $\calS^+_a$. The triality theory shows that if  $\inte \; \calS^+_a \neq \emptyset$,
     the global optimal solution can be obtained deterministically in polynomial time. Otherwise, the primal problem could be NP-hard
     and the nonlinear system could have chaotic behavior.
     A conjecture is proposed, which reveals the connection between chaos in nonlinear dynamics and NP-hardness in computer science.
     The method and the conjecture are verified by three examples. Computational results show that the conventional linear iteration methods
     can't produce reliable solutions due to the intrinsic numerical  error accumulations;
     while the canonical duality theory can be used not only to identify chaotic systems, but also for obtaining global optimal solutions to
     general nonlinear dynamic systems.

\bibliography{chaosbib}

\end{document}